\documentclass[12pt]{article}
\usepackage[latin1]{inputenc}
\usepackage{amssymb}
\usepackage{latexsym}
\usepackage{amscd}
\usepackage{longtable}
\usepackage{amsmath}
\usepackage{array}
\usepackage[all]{xy}
\usepackage{a4}
\usepackage{mathptmx}

\newtheorem{theorem}{Theorem}[section]
\newtheorem{cor}[theorem]{Corollary}
\newtheorem{lemma}[theorem]{Lemma}

\newtheorem{question}[theorem]{Question}
\newtheorem{proposition}[theorem]{Proposition}

\newtheorem{remark}{Remark}[section]

\newcommand{\Q}{\mathbb{Q}}
\newcommand{\Qbar}{\mathbb{Q}^{\operatorname{alg}}}
\newcommand{\Z}{\mathbb{Z}}

\newcommand{\C}{\mathbb{C}}
\newcommand{\R}{\mathbb{R}}

\newcommand{\Gal}{\operatorname{Gal}}

\newcommand{\SL}{\operatorname{SL}}

\newcommand{\End}{\operatorname{End}}

\newcommand{\Hom}{\operatorname{Hom}}
\newcommand{\Tr}{\operatorname{Tr}}
\newcommand{\Frob}{\operatorname{Frob}}
\newcommand{\Jac}{\operatorname{Jac}}
\newcommand{\Res}{\operatorname{Res}}
\newcommand{\Ver}{\operatorname{Ver}}

\newcommand{\id}{\operatorname{id}}
\newcommand{\im}{\operatorname{im}}
\newcommand{\ord}{\operatorname{ord}}
\newcommand{\pr}{\operatorname{pr}}

\newcommand{\Norm}{\operatorname{N}}
\newcommand{\gn}{\mathfrak n}
\newcommand{\gm}{\mathfrak m}

\newcommand{\ga}{\mathfrak a}

\newcommand{\gb}{\mathfrak b}

\newcommand{\gp}{\mathfrak p}

\newcommand{\gc}{\mathfrak c}
\newcommand{\gP}{\mathfrak P}

\newcommand{\cO}{{\mathcal O}}

\newcommand{\cE}{{\mathcal E}}
\newcommand{\cZ}{{\mathcal Z}}

\newcommand{\gmc}{{\overline{\mathfrak m}}}
\newcommand{\gpc}{{\overline{\mathfrak p}}}

\title{Modular elliptic directions \\
with Complex Multiplication \\
(with an application to Gross's elliptic curves)}

\author{Josep González and Joan-C. Lario}
\date{\today}

\begin{document}

\maketitle

\begin{abstract}
Let $A_f$ be the abelian variety attached by Shimura to a
normalized newform  $f\in S_2(\Gamma_1(N))$ and assume that $A_f$
has elliptic quotients. The paper deals with the determination of
the one dimensional subspaces (elliptic directions) in
$S_2(\Gamma_1(N))$ corresponding to the pullbacks of the regular
differentials of all elliptic quotients of $A_f$. For modular
elliptic curves over number fields without complex multiplication
(CM), the directions were studied by the authors in~\cite{gola01}.
The main goal of the present paper is to characterize the
directions corresponding to elliptic curves with CM. Then, we
apply the results obtained to the case $N=p^2$, for primes $p>3$
and $p\equiv 3 \bmod 4$. For this case we prove that if $f$ has
CM, then all optimal elliptic quotients of $A_f$ are also optimal
in the sense that its endomorphism ring is the maximal order of
$\Q(\sqrt{-p})$. Moreover, if $f$ has trivial Nebentypus then all
optimal quotients are Gross's elliptic curve $A(p)$ and its Galois
conjugates. Among all modular parametrizations $J_0(p^2)\to A(p)$,
we describe a canonical one and discuss some of its properties.
\end{abstract}

\section{Introduction}

Let $\Qbar$ be a fixed algebraic closure of $\Q$. An elliptic
curve $C$ defined over~$\Qbar$ is said to be modular if there is a
non-constant homomorphism $\pi\colon J_1(N) \to C$, where~$J_1(N)$
denotes the jacobian of the modular curve~$X_1(N)$.  Every modular
elliptic curve  over~$\Qbar$ is a quotient of some modular abelian
variety ${A_f}$ attached by Shimura to a normalized newform~$f$.
From now on, we shall always consider parametrizations $\pi\colon
J_1(N) \to C$ which factorize through such abelian varieties
$A_f$, called in this paper modular abelian varieties of {\it
elliptic type}.

A modular parametrization $\pi\colon J_1(N) \to C$ defined over a
number field $L\subseteq \Qbar$ induces an injection $\pi^* \colon
\Omega^1(C_{/L}) \hookrightarrow \Omega^1(J_1(N)_{/L})$. In what
follows, we shall identify $\Omega^1({J_1(N)}_{/L})$ with the
subspace of cusp forms in $S_2(\Gamma_1(N))$ whose $q$-expansion
lies in~$L[[q]]$, via $h\,dq/q\mapsto h$ where $q=\exp(2\pi i z)$.

The determination of the normalized cusp forms in
$S_2(\Gamma_1(N))$ associated with the pullbacks
$\pi^*(\Omega^1(C))$ was discussed by the authors in \cite{gola01}
for  elliptic curves without complex multiplication. In this
paper, we shall deal with the complex multiplication case that
needs techniques {\it ad hoc}. The present case is substantially
richer since it requires the intervention of class field theory as
well as the main theorem of complex multiplication.

Shimura shows in \cite{shimura71} that all elliptic curves with
complex multiplication (CM) are modular. Due to
Ribet~\cite{ribet77}, we know that $A_f$ has an elliptic quotient
with CM by an imaginary quadratic field $K\subset \Qbar$ if and
only if $f=f\otimes \chi$, where~$\chi$ is the quadratic Dirichlet
character attached to $K$. In this case, there is a primitive
Hecke character $\psi\colon I(\gm)\to \Qbar$ of conductor an ideal
$\gm$ of $K$ such that the $q$-expansion of the CM normalized
newform~$f$ is given by
$$
f = \sum_{(\ga,\gm)=1} \psi(\ga) q^{\Norm(\ga)}= \sum_{n=1}^\infty
a_n q^n \,.
$$
Here, $I(\gm)$ denotes the multiplicative group of fractional
ideals of $K$ relatively prime to~$\gm$, and the first summation
is over integral ideals. The level of $f$ is
$N=\Norm(\gm)\,|\Delta_K|$, the norm of $\gm$ times the absolute
value of the discriminant of $K$. We consider the number fields
$E_f=\Q(\{a_n\})$ and $E=\Q(\{\psi(\ga)\})$, generated by the
images of $\psi$. One has $E=E_f\cdot K$, and we shall denote by
$\Phi$ the set of its $K$-embeddings $E \hookrightarrow \Qbar$.
The number field $E$ is a CM field. Through the paper, for all CM
fields we shall denote by bar $\,\overline{\phantom{a}}\,$ the
canonical complex conjugation.

For future use, we recall that an abelian variety $Y$ is called
optimal quotient of an abelian variety $X$ over a field~$k$ if
there is a surjective morphism $\pi\colon X\rightarrow Y$ defined
over~$k$ whose kernel is an abelian variety. In this case, every
endomorphism of $X$ which leaves stable $\ker \pi$ induces an
endomorphism of $Y$. The property of being optimal quotient is
transitive. Hereafter, every $A_f$ is taken to be optimal quotient
of $J_1(N)$.

The plan of the paper is as follows. In section~2, we study the
decomposition of~$A_f$ over the quadratic field~$K$ for $f$ with
CM as before. This is an intermediate step necessary to determine
the elliptic directions we are interested in. We shall prove
\begin{theorem}\label{teoremaA} Let $f\in S_2(\Gamma_1(N))$ be a newform with
CM and keep the above notations. There is an abelian variety
$(A,\iota)$ of CM type $\Phi$ defined over $K$, with $\iota \colon
E \hookrightarrow \End_K^0(A)$, satisfying the following
properties:
\begin{itemize}
\item [(i)] $A$ is an optimal quotient of $A_f$ over $K$ and the pullback of $\Omega^1(A)$
corresponds with the subspace generated by $ \{ {}^{\sigma}f
 \colon \sigma\in\Phi \}$;
\item [(ii)]  $\iota(\psi(\ga))^*({}^{\sigma}f)={}^{\sigma}\psi(\ga)
{}^{\sigma}f$, for all $\ga\in I(\gm)$ and $\sigma\in \Phi$;
\item [(iii)] $\iota$ is an isomorphism;
\item [(iv)] if $\gp$ is a prime ideal of $K$ with $\gp \nmid N$, then the
lifting of the Frobenius endomorphism acting on the reduction of
$A$ mod $\gp$ is $\iota(\psi(\gp))$ or
$\iota(\overline{\psi(\overline{\gp})})$ depending on $K\nsubseteq
E_f$ or $K\subseteq E_f$, respectively.
\end{itemize}
\end{theorem}

We remark that the above abelian variety~$A$ is simple over $K$,
and that $A$ is $A_f$ over~$K$ when $K\nsubseteq E_f$, while~$A_f$
is isogenous over $K$ to $A\times \overline{A}$ when $K\subseteq
E_f$. To encode both cases of part~(iv) in Theorem~\ref{teoremaA},
we shall denote by $\psi'$ the primitive Hecke character
mod~$\gmc$ defined as
$$
\psi'(\ga)=
\begin{cases}
\psi(\ga)\,, & \text{ if $K \nsubseteq E_f$;}\\[8pt]
\overline{\psi(\overline{\ga})}\,, & \text{ if $K \subseteq E_f$.}
\end{cases}
$$
As it will be shown, one has $\gm=\gmc$ in the first case.

Then we study the splitting field of $A$; that is, the smallest
number field where all endomorphisms of $A$ are defined. We make
use of class field theory to build a certain abelian extension
$L/K$ attached to the Hecke character~$\psi'$; the field $L$ is a
cyclic extension of the Hilbert class field of~$K$ and it is
contained in the ray class field mod~$\overline{\gm}$. To simplify
notation, the Artin automorphism $(\frac{L/K}{\ga})$ in
$\Gal(L/K)$ will be often denoted by the same symbol representing
the ideal $\ga$. In particular, one has
$$
{}^\gp \beta \equiv \beta^{\Norm(\gp)} \pmod{\gP}
$$
for all $\beta\in {\mathcal{O}}_L$, where $\gP$ is an unramified
prime ideal of $L$ over a prime ideal $\gp$ of~$K$. The extension
$L/K$ is characterized by the property that $\ga$ viewed in
$\Gal(L/K)$ is trivial if and only if $\psi'(\ga)\in K^*$. The
main result of Section~3 is the following
\begin{theorem}\label{sistema}
Let $A$ be as above. Then,
\begin{itemize}
\item [(i)] there is an elliptic curve $C$ defined over $L$ with complex
multiplication by the ring of integers~${\mathcal O}_K$ and such
that $A$ is isogenous over $L$ to $C^{\dim A}$;
\item [(ii)]  the field $L$ is the smallest number
field satisfying $\End_{\Qbar}^0(A)=\End_{L}^0(A)$;
\item [(iii)] there is a one-cocycle $\lambda\colon I(\gmc) \to
L^*$ satisfying $\lambda(\ga) = \psi'(\ga)$ for all $\ga\in
I(\gmc)$ with $(\frac{L/K}{\ga})=\operatorname{id}$ in
$\Gal(L/K)$. The class of $\lambda$ in $H^1(I(\gmc), L^*)$ is
uniquely determined by this condition.
\end{itemize}
\end{theorem}

In view of (iii), the cohomology class of $\lambda$ depends
intrinsically on $A$, and we shall denote it by $[A]\in
H^1(I(\gmc), L^*)$. Section~4 is devoted to determining the
elliptic directions in $\Omega^1(A)$ in terms of $[A]$. To this
end, for each one-cocycle $\lambda\in [A]$ and $\sigma\in \Phi$,
we introduce the sums
$$
g_{\sigma}(\lambda):=\sum_{\ga\in \Gal(L/K)}
\frac{{}^{\ga^{-1}}\lambda(\ga)}{{}^\sigma \psi'(\ga)}\in
{}^\sigma E\cdot L\,,
$$
and also its $\Phi$-trace
$$
\operatorname{tr}_\Phi(\lambda) := \sum_{\sigma\in \Phi}
g_{\sigma}(\lambda) \in L \,.
$$

\pagebreak

\begin{theorem}\label{modular}
With the above notations, we have:

\begin{itemize}
\item [(1)] if $\,\sum_{n\geq 1} \gamma_n q^n \in S_2(\Gamma_1(N)) $ corresponds to
an elliptic direction attached to a modular parametrization
$\pi\in \Hom_L(A,C)$, then $\gamma_1\neq 0$.

\item [(2)] The following statements are equivalent:

\begin{itemize}
\item [(i)] the normalized cusp form $$
h= q+ \sum_{n\geq 2} \gamma_n q^n \in S_2(\Gamma_1(N))
$$
gives an elliptic direction attached to some $\pi\in \Hom_L(A,C)$;
\item [(ii)] there is a one-cocycle $\lambda \in [A]$ with $\operatorname{tr}_\Phi(\lambda)=[L\colon K]$
and such that
$$h=\frac{1}{[L:K]}\sum_{\sigma\in \Phi} g_\sigma(\lambda)\cdot {}^\sigma
f\,.$$
\end{itemize}
The $q$-expansion of this elliptic direction is then given by
$$h=
\left\{\begin{array}{ll}
\displaystyle{\sum_{(\ga,\gm)=1}{}^{\ga^{-1}}\lambda(\ga) \,q^{\Norm(\ga)}} &  \text{, if $K\nsubseteq E_f$;}\\[8 pt]
\displaystyle{\sum_{(\ga,\gm)=1}\frac{\Norm(\ga)}{\lambda(\overline
\ga)} \,q^{\Norm(\ga)}} & \text{, if $K\subseteq E_f$.}
 \end{array}\right.
$$

\noindent Moreover, all other elliptic directions are
$\iota(a)^*(h)$, for $a\in E^*$, and it holds the equality
$\iota(\psi'(\ga))^*h = {}^{\ga^{-1}}\lambda(\ga) {}^{\ga^{-1}} h$
for every $\ga\in I(\gmc)$.
\end{itemize}
\end{theorem}

We shall say that a one-cocycle $\lambda\in [A]$ is {\it modular}
if one has $\operatorname{tr}_\Phi(\lambda)=[L\colon K]$.
According to Theorem\ref{modular}, these are precisely the
one-cocycles that provide the elliptic directions. In Section~3,
we also describe how to obtain all modular one-cocycles in~$[A]$
explicitly way by means of a $K$-linear projector, and close the
section by raising some open questions.

\bigskip

In the last three sections, we deal with the particular case
concerning the level $N=p^2$ where $p>3$ is a prime with $p\equiv 3
\bmod 4$. The relevance of this case is in connection with the
elliptic curves $A(p)$ studied by Gross in \cite{gross80} and
\cite{gross82}. For convenience of the reader, we recall here its
definition. Let $K=\Q(\sqrt{-p})$ and let $\cO_K$ be its ring of
integers. Let $H$ denote the Hilbert class field of $K$, and let
$H_0=\Q(j(\cO_K))$ be its maximal real subfield. The elliptic curve
$A(p)$ is defined over $H_0$ and given by the Weierstrass equation
$$
y^2 = x^3 + \frac{mp}{2^4 \cdot 3}\,x - \frac{n p^2}{2^5 \cdot
3^3}\,,
$$
where $m$ and $n$ are the real numbers satisfying
$$
m^3 = j(\cO_K)\,,\quad n^2= \frac{j(\cO_K) -1728}{-p}\,,\quad
\operatorname{sgn}\, n = \left( \frac{2}{p} \right)\,.
$$
The elliptic curve $A(p)$ admits a global minimal model over~$H_0$
with discriminant $-p^3$ and whose invariants are $c_4= -mp$ and
$c_6=n p^2$.

Given any intermediate modular subgroup $\Gamma$ between
$\Gamma_1(p^2)$ and $\Gamma_0(p^2)$ and a normalized newform $f\in
S_2(\Gamma)$, we denote by $A_f^{(\Gamma)}$ its associated optimal
quotient of $\operatorname{Jac}(X_\Gamma)$, where~$X_\Gamma$
denotes the modular curve over $\Q$ attached to $\Gamma$.
According to this terminology, we have
$A_f^{(\Gamma_1(p^2))}=A_f$. In Section~5, we prove:

\begin{theorem}\label{Thm4} With the above notations, we have:
\begin{itemize}
\item[(i)]
for every positive divisor $d$ of $(p-1)/2$ there is a unique
abelian variety $A_f$ of~CM elliptic type in $J_1(p^2)$ such that
the Nebentypus of $f$ has order $d$; one has $K\not\subseteq E_f$,
$\dim A_f=[H:K]\varphi(d)$, where $\varphi$ is the Euler function,
and the splitting field of $A_f$ is the intermediate field between
$H$ and $H\cdot \Q(e^{2\pi\imath /p})$ of degree~$d$.

\item[(ii)] Let $f$ be a CM normalized newform in
$S_2(\Gamma_1(p^2))$ and let $\Gamma$ satisfying
$$\Gamma_1(p^2)\subseteq \Gamma \subseteq \Gamma_\varepsilon :=
\left\{ \left( \begin{array}{cc} a & b\\ c & d
\end{array}\right)\in \Gamma_0( p^2)\colon  \varepsilon (d)=1
\right\}\,,$$ where $\varepsilon$ is the Nebentypus of~$f$. Then,
all optimal elliptic quotients of $A_f^{(\Gamma)}$ have complex
multiplication by $\cO_K$. Moreover, if $f$ belongs to
$S_2(\Gamma_0(p^2))$, then all optimal quotients of
$A_f^{(\Gamma)}$ are defined over~$H$ and are precisely the
elliptic curve $A(p)$ and its Galois conjugates.
\end{itemize}
\end{theorem}

\bigskip

Among all modular parametrizations $J_0(p^2)\to A(p)$ one stands
out. In Section~6, we discuss this canonical parametrization and
give some of its arithmetical properties.

\begin{theorem}\label{Thm5}
Set $\gp=\sqrt{-p}\,\cO_K$.   Let $\delta\colon I(\gp) \to H$ be
the unique map defined by the conditions
$\delta(\ga)^{12}=\Delta(\cO)/\Delta(\ga)$ and
$\left(\frac{\Norm_{H/K}(\delta(\ga))}{\gp}\right)=1$. Let
$\omega$ denote a N\'eron differential of~$A(p)$. Then,

\begin{itemize}
\item[(i)] there is an optimal quotient $\pi\colon J_0(p^2)\longrightarrow A(p)$
such that $\pi^*(\omega)= c\, g(q)\,dq/q$ where the elliptic
direction is given by
$$g(q)= \sum_{(\ga,\gp)=1} \delta(\ga) q^{\Norm (\ga)}\in S_2(\Gamma_0(p^2))\,,$$
and $c\in \Z$ is a unit in $\Z[\frac{1}{2p}]$.
\item[(ii)] The complex lattice $\left\{ 2\pi\,i \int_\gamma g (z) dz\colon \gamma\in
H_1(X_0(p^2),\Z)\right\}$ is
$$
\frac{1}{c} \cdot i^{(p+1)/4}\,\, \cdot
\displaystyle{\sqrt[\leftroot{6}\uproot{32} h]{ \displaystyle{
{\rho\cdot(2\pi)^{(2 h+1-p)/4}\cdot\sqrt{p}^{(1-3h)/2}}}\cdot
\displaystyle{\prod_{\stackrel{1\leq
m<p}{\chi(m)=1}}\Gamma\left(\frac{m}{p}\right)}}}  \,  \cdot \cO_K
$$ where
$h$ is the class number of~$K$, the $h$-th root is taken to be
real, $\Gamma$ is the Gamma function, and $
\rho=\displaystyle{\prod_{\ga\in \Gal(H/K)} \frac{
\delta(\ga)}{\sqrt{\Norm(\ga)}} }$ is a positive unit of $H_0$.

\end{itemize}
\end{theorem}

Finally, in Section~7 we discuss how to compute the modular
elliptic directions for~$A_f$ when $f\in S_2(\Gamma_1(p^2))$ has
CM and its Nebentypus is nontrivial.

\section{The abelian variety $A$}
We shall adhere to the notations in the Introduction and prove
Theorem~\ref{teoremaA}. Let $\psi\colon I(\gm)\to \Qbar$ be the
fixed primitive Hecke character, and let
$$
f = \sum_{(\ga,\gm)=1} \psi(\ga) q^{\Norm(\ga)} =
\sum_{n=1}^{\infty} a_n q^n
$$
be its associated CM newform in $S_2(\Gamma_1(N))$. The optimal
quotient~$A_f$ of~$J_1(N)$ is defined over $\Q$ by $A_f = J_1(N) /
I_f( J_1(N))$, where $I_f(J_1(N)$ is the annihilator of $f$ in the
Hecke algebra acting on $J_1(N)$. In particular, the pullback
of~$\Omega^1({A_f}_{/\Qbar})$ is $\langle \{ {}^{\sigma} f\}
\rangle$ where~$\sigma$ runs over $\Gal(\Qbar/\Q)$. Recall that
$E_f=\Q(\{a_n\})$ and $E=\Q(\{ \psi(\ga) \})$. We fix an
isomorphism
$$
\iota\colon E_f \hookrightarrow \End_\Q^0(A_f)\,,
$$
in such a way that $\iota(a_n)$ corresponds to the Hecke operator
$T_n$ acting on $A_f$. The Nebentypus of $f$ is  the mod $N$
Dirichlet character $\varepsilon(d)=\chi(d)\psi((d))/d$, where
$\chi$ is the quadratic character attached to $K$. We recall that
$\iota(\varepsilon(d))$ is the diamond operator~$\langle d\rangle$
acting on $A_f$. One has
$$
\dim A_f = [E_f:\Q] =
\begin{cases}
\phantom{2\,}[E:K]\,, & \text{if $K\nsubseteq E_f$;}\\
2\,[E:K]\,, & \text{if $K\subseteq E_f$.}\\
\end{cases}
$$
Notice that $E=K\cdot E_f$. Now, we proceed to construct the
abelian variety $A$ over~$K$ of dimension $[E:K]$ with the
properties required in Theorem~\ref{teoremaA}. According to
Shimura's Proposition~8 in \cite{shimuraoptimal}, there exists
$u\in \End_K^0(A_f)$ such that
$$
u^*({}^\sigma f) = \sqrt{\Delta_K}\, \cdot {}^\sigma f
$$
for all $\sigma$ in $\Gal(\Qbar/\Q)$.  Here, the choice of the
square root $\sqrt{\Delta_K}$ fixes $u$ up to a sign. For the case
$K\nsubseteq E_f$, we let $A=A_f $ and extend $\iota$ to $E$,
$$
\iota\colon E \hookrightarrow \End_K^0(A_f)\,,
$$
via $\iota(\sqrt{\Delta_K})=u$. For the second case, we proceed as
follows. Since now $K\subseteq E_f$, there is $\alpha\in E_f$ such
that $\iota(\alpha)\in \End_\Q^0(A_f)$ acts as
$$
\iota(\alpha)^*({}^\sigma f)= {}^\sigma\sqrt{\Delta_K}\, \cdot
{}^\sigma f
$$
for all $\sigma$ in $\Gal(\Qbar/\Q)$. Then, consider the
involution $w:=\iota(\alpha) u^{-1}\in\End_K^0(A_f)$. Let $A$ be
the optimal quotient of $J_1(N)$ defined by $A_f/B$, where
$B=(1-w)A_f$. Clearly, the abelian variety $A$ is defined over
$K$, and $\Omega^1(A_{/K})$ is identified with $\langle
{}^{\sigma} f \rangle_{\sigma\in\Phi}$. Since $B$ is stable by
$\iota(E)$, the isomorphism $\iota\colon E \hookrightarrow
\End_\Q^0(A_f)$ induces in a natural way an embedding still
denoted by the same letter
$$
\iota\colon E \hookrightarrow \End_K^0(A)
$$
such that $\iota(\gamma)^*({}^\sigma f)={}^\sigma\gamma \cdot
{}^\sigma f$ for all $\gamma$ in~$E$ and all $K$-embeddings
$\sigma$ in~$\Phi$.  From the equality $\overline{w}=-w$, it
follows that $\overline{B}=(1+w)A_f$. Note that $\overline{B}$ is
$K$-isogenous to~$A$.

A case-by-case argument, employing that $\End^0_K(X)
\hookrightarrow \End^0_\Q(\Res_{K/\Q}(X))$ for any abelian variety
$X_{/K}$, shows that the abelian variety $A$ is $K$-simple in both
cases. Therefore, it follows that $\iota$ is an isomorphism. In
both cases, $A$ is an abelian variety of CM type $\Phi$ and
satisfies (i), (ii), and (iii) of Theorem~\ref{teoremaA}.

To conclude the proof, it remains to check the property~(iv)
relative to the Frobenius liftings. To this end, let $p$ be a
prime such that $p\nmid N$ and denote by $\Frob_p$ and ${\Ver}_p$
the Frobenius and the Verschiebung acting on the reduction of
$A_f$ modulo $p$, which satisfy $\Frob_p\cdot \Ver_p =p$. By the
Eichler-Shimura congruence, we know that
$$
\widetilde{T_p }= \Frob_p+ \Ver_p \cdot \widetilde{\langle
p\rangle}\,,
$$ where $\widetilde{T_p }$ and $\widetilde{\langle p\rangle}$ denote the reductions
of the Hecke operator $T_p$ and the diamond operator $\langle
p\rangle$ acting on $A_f$ mod $p$. Let us consider the two cases
separately.

Case $K\nsubseteq E_f$: first, assume that $p{\mathcal
O}_K=\gp\gpc$ splits in $K$. Since
$$ \iota(a_p)=\iota(\psi(\gp))+\iota(\psi (\overline{\gp}))\,,\quad
\iota( \psi (\gp))\cdot\iota(\psi (\overline{\gp}))= p \,\langle
p\rangle \,,$$ and $\widetilde{T_p }= \widetilde{\iota(a_p)}$, it
follows that the lifting of $\Frob_p$ is either $\iota(\psi(\gp))$
or $\iota(\psi(\gpc))$. Since a certain power of $\psi(\gp)$
belongs to $\gp$, one concludes that the lifting of
$\Frob_\gp=\Frob_p$ is $\iota(\psi(\gp))$. A similar argument
works when $p{\mathcal O}_K=\gp$ is inert in $K$, taking into
account that $\Frob_\gp=\Frob_p^2=-p\,\widetilde{\langle
p\rangle}=\widetilde{\iota(\psi((p))}$.

Case $K\subseteq E_f$: since $\iota(E)$ leaves the abelian
subvariety~$B$ stable, applying the same arguments as before, it
follows that~$\iota(\psi(\gp))$ is the lifting of~$\Frob_\gp$
acting on the reduction of~$B$ mod~$\gp$. Since~$A$ is
$K$-isogenous to~$\overline{B}$, the statement~(iv) holds in this
case as well. This completes the proof of Theorem~\ref{teoremaA}.

\vskip 0.3truecm

The following lemma will be used in the next sections.

\begin{lemma}\label{AfsobreK}
If  $K\not\subseteq E_f$, then $\gm=\gmc$.
\end{lemma}

\noindent {\bf Proof.} Since $K\not\subseteq E_f$, there is
$\sigma$ in $\Gal(\Qbar/K)$ such that ${}^\sigma f = \overline f$.
First, we prove that the Hecke characters ${}^\sigma\psi$ and
$\psi_c$ given by ${}^\sigma\psi(\ga)={}^\sigma(\psi(\ga))$ and
$\psi_c(\ga)=\overline{\psi(\overline{\ga})}$ coincide
on~$I(\gm\,\overline{\gm})$. Indeed, since ${}^\sigma
\varepsilon=\varepsilon^{-1}$ the assertion is immediate for prime
ideals $\gp\mid p$ when $p$ is inert. For the case that $p$ splits
completely in $K$, from the equalities ${}^\sigma
a_p=\overline{a_p}$ and
${}^\sigma\varepsilon(p)=\varepsilon^{-1}(p)$, that is,
$$ {}^\sigma \psi(\gp)+{}^\sigma
\psi(\overline\gp)=\psi_c(\gp)+\psi_c(\overline\gp)\quad
\text{and}\quad
{}^\sigma\psi(\gp)\cdot{}^\sigma\psi(\overline\gp)=\psi_c(\gp)\cdot\psi_c(\overline\gp)\,,
$$
it follows that ${}^\sigma\psi(\gp)$ is either $\psi_c(\gp)$ or
$\psi_c(\overline\gp)$. Again,  we obtain that ${}^\sigma
\psi(\gp)$ and $\psi_c(\gp)$ are equal because a certain power of
them lie in $\gp$. Both Hecke characters being primitive of
conductor $\gm$ and $\gmc$ respectively, we must have
$\gm=\overline{\gm}$.\hfill$\Box$

\section{Splitting field of $A$}
We first introduce an abelian extension $L/K$ that will play a key
role in the splitting of the abelian variety $A$ over $\Qbar$. Let
$\psi'$ be the primitive Hecke character mod~$\gmc$,
$$\psi'\colon I(\gmc)\to \Qbar\,,$$
given by $\psi'(\ga)=\psi(\ga)$ if  $K\nsubseteq E_f$ or
$\psi'(\ga)=\overline{\psi(\overline{\ga})}$ otherwise. We
consider the character $\eta\colon (\cO_K/\gmc)^*\to \Qbar$
defined by
$$
\eta(a)=\frac{\psi'((a))}{a}\,,\quad\text{for all $a\in \cO_K$
with $(a,\gmc)=1$}\,.
$$
One easily checks that $\eta$ is well-defined. Recall that the
existence of a Hecke character  mod~$\gmc$ is equivalent to the
condition that the composition $ \cO_K^*\hookrightarrow \cO_K
\rightarrow \cO_K/\gmc$ is a group monomorphism (see
\cite{shimura71}) and thus $\ker \eta\cap \cO_K^*=\{ 1\}$. By
class field theory, to the congruence subgroup
$$
P_{\eta}(\gmc)= \{ (a) \in I(\gmc) \colon \, a \bmod{\gmc} \in
\ker(\eta) \}
$$
it corresponds an abelian extension $L/K$. It is easy to check
that, for $\ga\in I(\gmc)$, one has $\ga\in P_{\eta}(\gmc)$ if and
only if $\psi'(\ga)\in K$. Let $K_{\gmc}$ denote the ray class
field of $K$ mod $\gmc$. Since the map $a\mapsto a\cO_K$ provides
an isomorphism between $\ker \eta $ and
$P_{\eta}(\gmc)/P_1(\gmc)$, by using the exact sequence
$$ 1\rightarrow \cO_K^{*} \rightarrow (\cO_K/\gmc)^{*} \rightarrow I(\gmc)/P_{1}(\gmc)
\rightarrow I(\cO_K)/P(\cO_K)\rightarrow 1\,,$$ one readily shows
that $L=K_{\gmc}^{\ker \eta}$ and $\Gal(L/H)$ is isomorphic to the
cyclic group $\im(\eta)/{\cO^*_K}$. Recall that here $H$ denotes
the Hilbert class field of~$K$ and, as usual, for any integral
ideal~$\gn$ we denote by $P(\gn)$ the subgroup of  $I(\gn)$ formed
by principal ideals and the subscript $1$ is for the subgroup of
principal ideals with a generator congruent to one mod~$\gn$. An
alternate route to define the extension $L/K$ is as follows. For
every $\sigma\in \Phi$, the character
$$
\chi_\sigma\colon \Gal(K_{\gmc}/K) \to {\Qbar}^*\,,\qquad
\chi_\sigma(\ga)=\frac{{}^\sigma \psi'(\ga)}{\psi'(\ga)}
$$
is well-defined via the Artin isomorphism $\Gal(K_{\gmc}/K)\simeq
I(\gmc)/P_1(\gmc)$. Due to the fact that $\bigcap_{\sigma\in \Phi}
{\ker\chi_\sigma}=P_{\eta}(\gmc)/P_1(\gmc)$, it follows that
$$L=K_{\gmc}^{\bigcap_{\sigma\in \Phi} {\ker\chi_\sigma}}\,.$$
Notice that $L/\Q$ is not necessarily a normal extension; in fact,
this is so if and only if $L=\overline{L}$.

\begin{proposition}\label{corba}
There is an elliptic curve $C$ defined over $L$ such that:
\begin{itemize}
\item [(i)] $\End_L(C)\simeq \cO_K$;
\item [(ii)] its Grossencharacter $\psi_C$ coincides with $\psi'\circ \Norm_{L/K}$;
\item [(iii)] $C$ is isogenous over $L$ to all its $\Gal(L/K)$-conjugates;
\item [(iv)] the abelian variety $A$ is isogenous over $L$ to the power $C^{[E:K]}$.
\end{itemize}
\end{proposition}

\noindent {\bf Proof.} The extreme cases~$L=H$ and~$L=K_{\gmc}$
are proved by Gross in~\cite{gross80} and  by de~Shalit
in~\cite{deshalit}, respectively. For the general case, one can
follow the same arguments. Let~$C_1$ be any elliptic curve
over~$L$ such that~$\End_L(C)\simeq \cO_K$. Let~$\gn$ be its
conductor. Once we fix an isomorphism $\theta\colon K \to
\End_L^0(C_1)$, we can consider the Grossencharacter
$\psi_{C_1}\colon I_L(\gn)\to K^*$ attached to the pair
$(C_1,\theta)$. For a prime ideal~$\gP$ of~$L$ relatively prime
to~$\gn$, we know that $\theta(\psi_{C_1}(\gP))$ is the lifting of
the $\gP$-Frobenius acting on the reduction of~$C_1$ mod~$\gP$.
Recall also that if~$\gP\in P_{1,L}(\gn)$ then
$\psi_{C_1}(\gP)=\Norm_{L/K}(\beta)$, where~$\gP=(\beta)$ with
$\beta\equiv 1\pmod{\gP}$.

By class field theory, the composition $\psi'\circ \Norm_{L/K}$
takes values in $K^*$ and the equality $\psi'\circ
\Norm_{L/K}(\gP)=\Norm_{L/K}(\beta)$ holds for $\gP=(\beta)$ with
$\beta\equiv
 1\pmod{\gmc\cO_L}$. Therefore, the quotient
$(\psi'\circ \Norm_{L/K})/\psi_{C_1}$ defines a character
$\delta\colon I_L(\gn\gmc\cO_L)/P_{1,L}(\gn\gmc\cO_L) \to \cO_K^*$
of finite order. The twist $C:=C_1\otimes \delta$ satisfies~(i)
and~(ii). Now, (iii) follows from the fact that
$\psi_{C}=\psi_{{}^\ga C}$ for all $\ga\in \Gal(L/K)$ due to (ii).

Now, we check (iv). By Faltings's criterion (for instance, see \S
2 Corollary~2 of~\cite{cornell-silverman}), it suffices to prove
that for every prime~$\gP$ of~$L$ not dividing~$N$ nor the
conductor of~$C$, the reductions of the abelian varieties~$A$ and
$C^{\dim A}$ modulo~$\gP$ are isogenous over the residue field
$\cO_L/\gP$. We write~$\gp^f=\Norm_{L/K} \gP$, where with no risk
of confusion now~$f$ is the residue degree of~$\gP$ over~$K$. On
the one hand, the characteristic polynomial of the
endomorphism~$\Frob_{\gP}$ acting on the $l$-adic Tate module of
the reduction of~$A/L$ modulo~$\gP$, for a prime~$l\neq p$, is the
characteristic polynomial of the complex representation
of~$\iota(\psi'(\gp^f))$:
$$
 P_{A,\gP}(x)=
\prod_{\sigma\in
\Phi}(x-{}^{\sigma}\psi'(\gp^f))(x-\overline{{}^{\sigma}\psi'(\gp^f)})\,.
$$
On the other hand, the corresponding Frobenius characteristic
polynomial for $C$ at~$\gP$ is
$$
P_{C,\gP}(x)=(x-\psi_{C}(\gP))(x-\overline{\psi_{C}(\gP)}\,)=
(x-\psi'(\gp^f))(x-\overline{\psi'(\gp^f)})\,.
$$
Since  $\psi'(\gp^f)$ belongs to $K$,  we obtain
$P_{A,\gP}(x)=P_{C,\gP}(x)^{\dim A}$. Thus, $A$ is isogenous
over~$L$ to $C^{\dim A}$. \hfill$\Box$

\begin{proposition}
The field $L$ is the smallest number field satisfying
$\End_{\Qbar}^0(A)=\End_{L}^0(A)$.
\end{proposition}

\noindent {\bf Proof.} Since $A$ is isogenous over $L$ to the
$[E:K]$-th power of the elliptic curve~$C$, we have
$\End_{\Qbar}^0(A)=\End_{L}^0(A)$. That $L$ is the smallest number
field with this property can be deduced from the following fact.
For every $\varphi\in \End_L^0(A)$, one has the explicit version
of the Skolem-Noether theorem:
$$
{}^\gp \varphi  = \iota(\psi'(\gp)) \cdot \varphi \cdot
\iota(\psi'(\gp))^{-1}\,,
$$
for all $\gp\in I(\gmc)$ not dividing~$N$. To check this equality,
it is enough to verify that it holds reduced modulo a prime ideal
$\gP$ of $L$ over $\gp$. The smallest field of definition for all
endomorphisms of $A$ is the fixed field $L^G$, where
$$G=\{\nu \in \Gal (L/K):
{}^{\nu}\phi=\phi \quad\text{for all }\phi\in \End_{L}^0(A)\}\,.
$$
By the \v Cebotarev density theorem, every $\nu$ in $\Gal(L/K)$
can be written as $\nu =(\frac{L/K}{\gp} )$ for some prime ideal
$\gp$ relatively prime to~$N$. We have that $\nu \in G$ if and
only if $\iota(\psi'(\gp))$ is in the center of $\End_{L}^0(A)$;
that is, when  $\psi' (\gp)\in K$ and this fact implies that $\gp$
splits completely in $L$, so that $\nu=\id$.\hfill$\Box$

\vskip 0.2truecm

Let $C$ be an elliptic curve defined over $L$ as in
Proposition~\ref{corba}. The main theorem of complex
multiplication (Theorem 5.4 in \cite{shimurabook}) implies the
existence of a system of isogenies $\{\mu_\ga\colon C \to {}^\ga
C\}$ over $L$, $(\ga,\gmc)=1$, satisfying the following
properties:
\begin{itemize}
\item [(i)] $\mu_{\ga\gb}={}^\ga \mu_\gb\, \mu_\ga$;
\item [(ii)] if  $C$ has good reduction at a prime ideal $\gP\mid \gp$,
then $\mu_\gp$ is the lifting of the Frobenius map between the
reductions of $C$ and ${}^\gp C$ mod $\gP$.
\end{itemize}
Attached to the system of isogenies $\{\mu_\ga\}$, a one-cocycle
can be defined as follows (see also~\cite{go-ca}). For a non-zero
regular differential $\omega$ in $\Omega^1(C_{/L})$, let
$\lambda_\omega\colon I(\gmc)\to L^*$ be the map given by
$$
\mu_\ga^*({}^\ga \omega)=\lambda_\omega(\ga) \omega\,,
$$
where ${}^\ga \omega$ denotes the differential in ${}^\ga C$
corresponding to $\omega$ by conjugation. It follows that
$\lambda_\omega$ is a one-cocycle, and for all $u\in L^*$ one has
$$\lambda_{u\omega} (\ga)=\lambda_\omega(\ga){}^\ga u/u\,.$$
Clearly, the class of $\lambda_\omega$ in $H^1(I(\gmc),L^*)$ does
not depend on the particular choice of~$\omega$. Note that if
$\ga\in P_\eta(\gmc)$, then we have
$\lambda_\omega(\ga)=\psi'(\ga)$. The class $\lambda_\omega$ in
$H^1(I(\gmc),L^*)$ can be characterized from $\psi'$ as follows:

\begin{proposition}\label{lambda}
Let $\lambda \colon I(\gmc)\to L^*$ be any one-cocycle satisfying
$\lambda(\ga) = \psi'(\ga)$ for all $\ga\in I(\gmc)$ with
$(\frac{L/K}{\ga})=\operatorname{id}$ in $\Gal(L/K)$. Then,
$[\lambda]=[\lambda_\omega]$.
\end{proposition}

\noindent {\bf Proof.} Assume that $\lambda\in H^1(I(\gmc),L^*)$
satisfies $\lambda (\ga)=\psi'(\ga)$ for all $\ga\in P_\eta
(\gmc)$. The quotient $\lambda/\lambda_\omega$ defines a
one-cocycle in $H^1(\Gal(L/K),L^*)$. By Hilbert's 90 theorem, we
know that there is $u\in L^*$ such that
$\lambda(\ga)/\lambda_\omega(\ga)={}^\ga u/u$ for all $\ga\in
I(\gmc)$. Thus, we have $[\lambda]=[\lambda_\omega]$. \hfill$\Box$

\vskip 0.3truecm

This completes the proof of Theorem~\ref{sistema} in the
Introduction. From now on, we shall denote by~$[A]$ in
$H^1(I(\gmc),L^*)$ the cohomology class of~$\lambda_\omega$.

\section{Modular one-cocycles and elliptic directions}

In this section we keep the notations as above and tackle the
problem to determine the elliptic directions in $\Omega^1(A)$. The
goal is to prove Theorem~\ref{modular} that will be deduced from
the next three Propositions after the following

\begin{lemma}\label{gamma1}
Let $\pi\in \Hom_L(A,C)$ be a non-constant modular
parametrization, and let $\omega\in \Omega^1(C_{/L})$ be any
non-zero regular differential. Denote by
$$
h= \sum_{n\geq 1} \gamma_n q^n \in S_2(\Gamma_1(N))
$$
the cusp form associated with the pullback $\pi^*(\omega)$. Then,
\begin{itemize}
\item [(i)] $\gamma_1\in L^*$;
\item [(ii)] for all $\ga\in I(\gmc)$ relatively prime to $N$, one has $\iota(\psi'(\ga))^*h = {}^{\ga^{-1}}\lambda_\omega(\ga) {}^{\ga^{-1}}
h$;\\[1pt]
\item [(iii)] we have the identity
$ \displaystyle{ h =\frac{1}{[L:K]} \sum_{\ga\in\Gal(L/K)}
\sum_{\sigma\in \Phi}
\frac{{}^{\ga^{-1}}\lambda_\omega(\ga)}{{}^\sigma \psi'(\ga)}
{}^{\ga^{-1}} h\,. }$; \\[1pt]
\item [(iv)] $\{\psi'(\ga_i)\}$ is a $K$-basis of $E$ if and only if
$\{ {}^{\ga_i^{-1}} h\}$ is an $L$-basis of $\Omega^1(A_{/L})$.
\end{itemize}
\end{lemma}

\noindent {\bf Proof.} (i) Since $\pi$ and $\omega$ are defined
over $L$, the cusp form $h$ associated with $\pi^*(\omega)$ has
$q$-expansion $\sum_{n\geq 1} \gamma_n q^n$ with coefficients
in~$L$. Since the abelian variety $A$ is simple over $K$, we have
that $A$ is a $K$-factor of the Weil restriction $\Res_{L/K}(C)$.
Thus, the set $\{ {}^\ga h \colon \ga\in \Gal(L/K)\}$ generates
$\Omega^1(A_{/L})$. This implies $\gamma_1\neq 0$.

(ii) It is enough to consider the case when $\ga=\gp$ is a prime
ideal not dividing~$N$. Then, the claim follows from the
commutativity of the diagram
$$ \xymatrix{
A \ar[rr]^-{\iota(\psi'(\ga))} \ar[dd]_{{}^{\ga^{-1}}\pi}   && A
\ar[dd]_{\scriptstyle{\pi}}
\\  \\
{}^{\ga^{-1}}C \ar[rr]^-{{}^{\ga^{-1}}\mu_\ga}   &&  C}
\label{commute}
$$
due to the fact that $\iota(\psi'(\gp))$ and
${}^{\gp^{-1}}\mu_\gp$ are liftings of the corresponding
$\gp$-Frobenius morphisms at a prime ideal $\gP\mid\gp$ of $L$.

(iii) Write $h=\sum_{\nu\in \Phi} c_\nu {}^{\nu}f$, with $c_\nu\in
\Qbar$. By applying (ii), for all $\sigma\in\Phi$ and
$\ga\in\Gal(L/K)$, one has
$$
\frac{{}^{\ga^{-1}}\lambda_\omega
(\ga)}{{}^{\sigma}\psi'(\ga)}\,{}^{\ga^{-1}}h=\frac{1}{{}^\sigma
\psi' (\ga)}\left ( \sum_{\nu\in\Phi} c_\nu
{}^\nu\psi'(\ga)\,{}^\nu f \right)=\sum_{\nu\in\Phi} c_\nu (
\chi_\nu \cdot {\chi_\sigma }^{-1})(\ga)\,{}^\nu f\,.
$$
Thus, it holds
$$
\sum_{\ga} \sum_{\sigma} \frac{{}^{\ga^{-1}}\lambda_\omega
(\ga)}{{}^{\sigma}\psi'(\ga)}\,{}^{\ga^{-1}}h=\sum_{\sigma,\nu}\sum_\ga
c_\nu ( \chi_\nu \cdot {\chi_\sigma }^{-1})(\ga)\,{}^\nu f=[L:K]\,
\sum_\nu c_\nu {}^\nu f=[L:K]\, h\,.
$$

(iv) If $\{\psi'(\ga_1),\dots, \psi'(\ga_r)\}$ is a $K$-basis of
$E$, then for every $\ga\in I(\gmc)$ we can write
$\psi'(\ga)=\sum_{i=1}^r \alpha_i \psi'(\ga_i)$ with $\alpha_i\in
K$. Thus, we obtain
$$
{}^{\ga^{-1}} \lambda_\omega(\ga) {}^{\ga^{-1}} h =
\iota(\psi'(\ga))^*(h)=\sum_{i=1}^r \alpha_i \iota(\psi'(\ga_i))^*
h = \sum_{i=1}^r \alpha_i {}^{\ga_i^{-1}} \lambda_\omega(\ga_i)
{}^{\ga_i^{-1}}h\,.
$$
Since $\{ {}^\ga h \colon \ga\in \Gal(L/K)\}$ generates
$\Omega^1(A_{/L})$ and $\dim(A)=[E:K]$, it follows that $\{
{}^{\ga_1^{-1}} h, \dots, {}^{\ga_r^{-1}} h \}$ is a $L$-basis of
$\Omega^1(A_{/L})$.

Conversely, assume that $\{ {}^{\ga_1^{-1}} h ,\dots,
{}^{\ga_r^{-1}} h \}$ is a $L$-basis of $\Omega^1(A_{/L})$. By
using part~(ii), if $\sum_{i=1}^r \alpha_i \psi'(\ga_i) = 0$ for
some $\alpha_i\in K$, then $\sum_{i=1}^r \alpha_i \,
{}^{\ga_i^{-1}} \, \lambda_\omega(\ga_i){}^{\ga_i^{-1}}h = 0$.
This implies that all $\alpha_i=0$. Since $\dim(A)=[E:K]=r$, the
proof is done. \hfill$\Box$

\vskip 0.2truecm

Due to  part (i) in the above Lemma~\ref{gamma1}, there is a
unique $\omega\in \Omega^1(C_{/L})$ such that the pullback
$\pi^*(\omega)$ gives a normalized cusp form, say
$$
h= q+\sum_{n\geq 2} \gamma_n q^n\,.
$$
The particular $\lambda_\omega$ will be called {\it modular with
respect to} $\pi$ or, simply, $\pi$-{\it modular}.  For every
$1$-cocycle $\lambda\in [A]$, we consider the following sums. Let
$\sigma\in \Phi$, and set
$$
g_\sigma(\lambda):=\sum_{\ga\in \Gal(L/K)}
\frac{{}^{\ga^{-1}}\lambda(\ga)}{{}^\sigma \psi'(\ga)}\,.
$$
Notice that $g_\sigma(\lambda)$ is well-defined and
$g_\sigma(\lambda) \in {}^\sigma E \cdot L$.

\begin{remark}\label{global}
The sum $g_\sigma(\lambda)$ can be interpreted as a sort of Gauss
sum, in the sense that we have
$$ g_\sigma(\lambda)= \sum_{\ga\in \Gal(L/K)} \chi_\sigma^{-1}(\ga) u_\ga$$
where $u_\ga =  {}^{\ga^{-1}}\lambda(\ga)/ \psi'(\ga)$. If $C$
admits a global minimal Weierstrass equation over~$L$, then the
one-cocycle $\lambda$ attached to a N\'eron differential satisfies
the capitulation property $\lambda(\ga)\cO_L = \ga\cO_L$ (see
remark~10.3 in \cite{go-ca}). Then, $u_\ga^e$ is an unit in
$\cO_L^*$ where $e$ is the order of~$\ga$ in~$\Gal(L/K)$.
\end{remark}

We shall denote the $\Phi$-trace of $g_\sigma(\lambda)$ by
$$\operatorname{tr}_\Phi(\lambda)=\sum_{\sigma\in\Phi} g_\sigma(\lambda) \in L\,.$$

\begin{remark}
Recall that we have defined $\lambda\in [A]$ to be modular if
$\operatorname{tr}_\Phi(\lambda)=[L\colon K]$ in the Introduction.
As it will be shown, both terms (modular and $\pi$-modular) turn
out to be equivalent.
\end{remark}

For every $\gamma\in L^*$ and $\lambda\in [A]$, let
$\lambda_\gamma$ denote the twisted one-cocycle in $[A]$ given by
$\lambda_\gamma(\ga)=\lambda(\ga) \gamma/{}^\ga \gamma$. Writing
$\lambda=\lambda_\omega$ with some $\omega\in \Omega^1(C_{/L})$,
then $\lambda_{\gamma}=\lambda_{\frac{1}{\gamma}\omega}$. We shall
need the following lemma.

\begin{lemma}\label{gauss}
For all $\ga\in I(\gmc)$ and $\sigma\in \Phi$, one has
\begin{itemize}
\item [(i)] $\displaystyle{g_{\sigma}(\lambda_{{}^{\ga^{-1}}\lambda(\ga)})\cdot  {}^{\ga^{-1}}\lambda(\ga)=
g_\sigma(\lambda)\cdot {}^\sigma \psi'(\ga)}$;
\item [(ii)] $\operatorname{tr}_\Phi (\lambda_{ {}^{\ga^{-1}} \lambda(\ga)}) =
{}^{\ga^{-1}}\operatorname{tr}_\Phi (\lambda)$.
\end{itemize}
\end{lemma}

\noindent {\bf Proof.} It follows straightforward from the
definitions and by using the cocycle relations
for~$\lambda$.\hfill$\Box$

\begin{proposition}\label{p1}
Assume that $\lambda\in [A]$ is modular with respect to $\pi\in
\Hom_L(A,C)$. Then, $\operatorname{tr}_\Phi(\lambda)=[L:K]$ and
$$
h=\frac{1}{[L:K]}\sum_{\sigma\in\Phi} g_{\sigma}(\lambda)\cdot
{}^\sigma f
$$
is the normalized elliptic direction in $\pi^*(\Omega^1(C_{/L}))$.
\end{proposition}

\noindent {\bf Proof.} Since $\lambda$ is $\pi$-modular, there is
a non-zero regular differential $\omega\in \Omega^1(C_{/L})$ such
that $\pi^*(\omega)$ is a normalized cusp form $h=q+\sum_{n\geq 2}
\gamma_n q^n$ and $\lambda=\lambda_\omega$. By comparing the first
Fourier coefficient in the equality at Lemma~\ref{gamma1}~(iii),
we have that $\operatorname{tr}_\Phi(\lambda)=[L:K]$. For every
$\sigma\in\Phi$, set
$$
F_\sigma=\sum_{\gb\in\Gal(L/K)}\frac{{}^{\gb^{-1}}\lambda(\gb)}{{}^{\sigma}\psi'
(\gb)} \, {}^{\gb^{-1}}h\,.
$$
Also by Lemma~\ref{gamma1}~(iii), we know that $\sum_{\sigma
\in\Phi}F_\sigma=[L:K]\,h$. From the equality
$\iota(\psi'(\ga))^*(\, {}^{\gb^{-1}} h)={}^{(\gb\cdot
\ga)^{-1}}\lambda(\ga)\,{}^{(\gb\cdot\ga)^{-1}} h$, one obtains
$$
\begin{array}{l@{\,=\,}l}
\iota(\psi'(\ga))^* (F_\sigma)&
\displaystyle{\sum_{\gb\in\Gal(L/K)}\frac{{}^{\gb^{-1}}\lambda(\gb)}{{}^{\sigma}\psi'
(\gb)}\, {}^{(\gb\cdot \ga)^{-1}}\lambda(\ga)\,{}^{(\gb\cdot\ga)^{-1}} h} \\[10pt]
&
\displaystyle{\sum_{\gb\in\Gal(L/K)}\frac{{}^{(\gb\cdot\ga)^{-1}}\lambda(\gb\cdot
\ga)}{{}^{\sigma}\psi' (\gb)}\, {}^{(\gb\cdot \ga)^{-1}} h
={}^\sigma \psi'(\ga)\,F_\sigma\,.}
\end{array}
$$
Therefore, $F_\sigma$ and ${}^\sigma f$ differ by a scalar
multiple. Since the $q$-expansion of $F_\sigma$ begins as
$g_\sigma(\lambda)\,q+\cdots\,,$ it follows that
$F_\sigma=g_\sigma(\lambda)\cdot {}^\sigma f$, and then
$\displaystyle{h=\frac{1}{[L:K]}\sum_{\sigma\in \Phi}
g_\sigma(\lambda) \cdot {}^\sigma f}$.~\hfill$\Box$

\vskip 0.3truecm

Now, we shall prove that the modular one-cocycles $\lambda$ in
$[A]$ with respect to some modular parametrization $\pi$ are
precisely those that satisfy the trace
condition~$\operatorname{tr}_\Phi(\lambda)=[L:K]$. To this end,
for a given one-cocycle $\lambda\in [A]$ (not necessarily
modular), let us consider the $K$-linear map
$\operatorname{pr}\colon L \to L$,
$$
\operatorname{pr}(u):= \sum_{\ga\in\Gal(L/K)} \left(
\sum_{\sigma\in\Phi} \frac{1}{{}^\sigma \psi'(\ga)}\right)
{}^{\ga^{-1}}\lambda(\ga) {}^{\ga^{-1}} u=
\begin{cases}
u \cdot \operatorname{tr}_\Phi(\lambda_u) & \text{, if $u\neq
0$;}\\[6pt]
0 & \text{, otherwise.}
\end{cases}
$$

Consider the eigenspace ${\mathcal M} = \{ u\in L\colon
\pr(u)=[L:K]\cdot u \}$. Notice that $\lambda_u$ is modular if and
only if~$u\in {\mathcal M}\backslash \{ 0\}$. In particular, we
know that $\dim_K( {\mathcal M} )>0$ and it does not depend on the
particular choice of $\lambda\in [A]$ used to define the
$K$-linear map~$\operatorname{pr}$.

\begin{proposition}\label{projector}
One has
\begin{itemize}
\item [(i)] $\operatorname{pr}^2 = [L:K]\, \operatorname{pr}$;
\item [(ii)] $\dim_K({\mathcal M})= [E:K]$;
\item [(iii)] if $\lambda$ is modular,
then ${\mathcal M}=\langle \{ {}^{\ga^{-1}}\lambda(\ga)\}
\rangle_K$ where $\ga$ runs over $\Gal(L/K)$.
\end{itemize}
\end{proposition}

\noindent {\bf Proof.} The first claim comes from the computation:
$$
\begin{array}{l@{\,=\,}l}
\operatorname{pr}^2(u) & \displaystyle{\sum_\ga \left( \sum_\sigma
\frac{1}{{}^\sigma\psi'(\ga)}\right) {}^{\ga^{-1}} \lambda(\ga)
\!\!\!\!{\phantom{\int}}^{\ga^{-1}} \left[ \sum_\gb \left(
\sum_\tau \frac{1}{{}^\tau\psi'(\gb)}\right) {}^{\gb^{-1}}
\lambda(\gb) {}^{\gb^{-1}} u
\right]} \\[10pt]
& \displaystyle{\sum_\ga \sum_\gb \left( \sum_\sigma
\frac{1}{{}^\sigma\psi'(\ga)}\right)
\left( \sum_\tau \frac{1}{{}^\tau\psi'(\gb)}\right) {}^{(\ga\gb)^{-1}} \lambda(\ga\gb) {}^{(\ga\gb)^{-1}}u} \\[10pt]
& \displaystyle{\sum_\ga \sum_\gb \left( \sum_\sigma
\frac{1}{{}^\sigma\psi'(\ga)}\right)
\left( \sum_\tau \frac{1}{{}^\tau\psi'(\ga^{-1}\gb)}\right) {}^{\gb^{-1}} \lambda(\gb) {}^{\gb^{-1}}u} \\[10pt]
& \displaystyle{\sum_\gb \sum_\ga  \left( \sum_{\sigma,\tau}
(\chi_\sigma\chi_{\tau}^{-1})(\ga)\right)
 \frac{   {}^{\gb^{-1}} \lambda(\gb)  }{  {}^\tau\psi'(\gb)}  {}^{\gb^{-1}} u} \\[14pt]
& [L:K]\,\operatorname{pr}(u)\,.
\end{array}
$$

Let us  prove (ii) and (iii) simultaneously. Since $\dim
_K({\mathcal M})$ is independent of the one-cocycle $\lambda$
chosen in $[A]$, we can (and do) assume that $\lambda$ is modular.
Set
$$
h=\frac{1}{[L:K]}\sum_{\sigma\in\Phi} g_{\sigma}(\lambda)\cdot
{}^\sigma f=1+ \sum_{n>1}\gamma_n \,q^n\,.
$$
Let $W= \langle  \{ {}^{\ga^{-1}}\lambda(\ga) \} \rangle_K$ where
$\ga$ runs over $\Gal(L/K)$. We need to show that $W={\mathcal M}$
and $\dim_K(W)=[E:K]$.
 Choose
 $\ga_1,\dots,\ga_r\in I(\gmc)$ such that
 $\{\psi'(\ga_1),\dots, \psi'(\ga_r) \}$
 is a $K$-basis of $E$. We claim that
 $\{ {}^{\ga_1^{-1}}\lambda (\ga_1) , \dots, {}^{\ga_r^{-1}}\lambda (\ga_r) \}$
 is a $K$-basis of $W$.
 Indeed, if
 $\sum_{i=1}^r \alpha_i \,{}^{\ga_i^{-1}}\lambda (\ga_i)=0$ for some
 $\alpha_i$ in~$K$, then consider $\alpha:=\sum_{i=1}^r\alpha_i
\psi'(\ga_i)\in E$. It is easy to check that $\iota(\alpha)^*(h)
=\sum_{n\geq 1}\gamma_n' \,q^n$ with $\gamma_1'=0$.  This forces
$\alpha= 0$, since otherwise we get a contradiction from
Lemma~\ref{gamma1}~(i) applied to $\iota( \alpha)^*(h)$.
Therefore, all $\alpha_i=0$ which implies that $
{}^{\ga_1^{-1}}\lambda (\ga_1) , \dots, {}^{\ga_r^{-1}}\lambda
(\ga_r)$ are linearly independent. Now, for every ideal  $\ga\in
I(\gmc)$, one has $\psi'(\ga)=\sum_{i=1}^r\alpha_i\psi'(\ga_i)$
for some $\alpha_i\in K$. By taking $q$-expansions in the equality
 $$
{}^{\ga^{-1}}\lambda(\ga){}^{\ga^{-1}}h=\sum_{i=1}^r
\alpha_i\,{}^{\ga_i^{-1}} \lambda(\ga_i) \,{}^{\ga_i^{-1}} h\,,
 $$
we obtain ${}^{\ga^{-1}}\lambda(\ga)=\sum_{i=1}^r
\alpha_i\,{}^{\ga_i^{-1}}\lambda(\ga_i)$. So far, we have $\dim_K
(W) =[E:K]$ and the inclusion $W \subseteq {\mathcal M}$ follows
from Lemma~\ref{gauss}~(ii).

To easy notation, set $u_i={}^{\ga_i^{-1}}\lambda(\ga_i)$ for
$1\leq i \leq r$ and let us show that they generate~${\mathcal
M}$. For any nonzero $u\in{\mathcal M}$, consider the normalized
cusp form
$$
h_u=\frac{1}{[L:K]}\sum_{\sigma\in\Phi} g_{\sigma}(\lambda_u)\cdot
{}^\sigma f \,.
$$
Since $\{h_{u_1},\cdots,h_{u_n} \}$ is a $L$-basis of
$\Omega^1(A_{/L})$ by Lemma~\ref{gamma1}~(iv), there are
$\gamma_i\in L$ such that $h_u=\sum_{i=1}^r \gamma_i \, h_{u_i}$.
Notice that $\sum_{i=1}^r \gamma_i=1$. By applying $\iota
(\psi'(\ga))^*$ to $h_u$, and then conjugate by $\ga$,  we obtain
$$
\lambda_{u}(\ga) h_u=\sum_{i=1}^r {}^\ga \gamma_i\,
\lambda_{u_i}(\ga)\, h_{u_i}\,.
$$
Therefore, we have
$$\gamma_i={}^\ga \gamma_i
\frac{\lambda_{u_i}(\ga)}{\lambda_u(\ga)} = {}^\ga \gamma_i
\frac{{}^\ga u}{{}^\ga u_i} \frac{u_i}{u} $$ for all $\ga$ and
$1\leq i \leq r$. That is, $\beta_i:=\gamma_i\, u/u_i\in K$. Then,
$u=\sum_{i=1}^r \beta_i \,u_i$ since $\sum_{i=1}^r \gamma_i=1$.
The statement (iii) follows.\hfill$\Box$

\begin{proposition}\label{p3}
Let $\lambda'\in [A]$ such that $\operatorname{tr}_\Phi (\lambda')
= [L: K]$. Then, $\lambda'$ is modular with respect to
some~$\pi'\in \Hom_L(A,C)$.
\end{proposition}

\noindent {\bf Proof.}  We shall prove that there is $\pi'\in
\Hom_L(A,C)$ and $\omega'\in \Omega^1(C_{/L})$ such that
$\pi'^*(\omega')$ corresponds to the normalized cusp form
$$
h'=\frac{1}{[L:K]} \sum_{\sigma\in \Phi} \sum_{\ga\in \Gal(L/K)}
\frac{{}^{\ga^{-1}}\lambda'(\ga)}{{}^\sigma \psi'(\ga)} \cdot
{}^\sigma f\,.
$$
Consider any non-constant $\pi\in \Hom_L(A,C)$ and take $\omega\in
\Omega^1(C_{/L})$ such that $\pi^*(\omega)$ corresponds to the
normalized cusp form
$$
h=\frac{1}{[L:K]} \sum_{\sigma\in \Phi} g_\sigma(\lambda)\cdot
{}^\sigma f\,,
$$
where $\lambda=\lambda_\omega$. Let
$L=\ker(\operatorname{pr})\oplus {\mathcal M}$ be the
decomposition corresponding to the projector $\operatorname{pr}$
attached to $\lambda$. Now, there is $\gamma\in {\mathcal M}$ such
that $\lambda'=\lambda_\gamma$ and
$$
h'=\frac{1}{[L:K]} \sum_{\sigma\in \Phi}
g_{\sigma}(\lambda_\gamma) \cdot {}^\sigma f
$$
with $\gamma =  \sum_{\ga\in \Gal(L/K)} r_\ga {}^{\ga^{-1}}
\lambda(\ga)$ for some $r_\ga\in K$ due to
Proposition~\ref{projector}~(iii). We claim that
\begin{equation}\label{igualtat}
\left( \sum_{\ga\in \Gal(L/K)} r_\ga {}^{\ga^{-1}}
\lambda(\ga)\right) h' = \iota \left( \sum_{\ga\in \Gal(L/K)}
r_\ga  \psi'(\ga)\right)^* h\,.
\end{equation}
Letting $\Psi= \iota\left(\sum_{\ga\in \Gal(L/K)} r_\ga
\psi'(\ga)\right) \in \End_K^0(A)$, then it follows
$$
h' = \Psi^*\left( \pi^*(\frac{1}{\gamma} \omega) \right) = (\pi
\circ \Psi)^* \left(\frac{1}{\gamma} \omega \right)\,,
$$
which implies that $\lambda'$ is modular. To check
(\ref{igualtat}), we use Lemma~\ref{gauss}~(i):
$$
\begin{array}{ll}
\gamma \, h' &  =\displaystyle{ \frac{1}{[L:K]}\sum_{\sigma}
\sum_{\gb}
 \frac{{}^{\gb^{-1}}\lambda(\gb) {}^{\gb^{-1}}\gamma  }{{}^\sigma \psi'(\gb)} \,{}^\sigma f} \\[10pt]
& = \displaystyle{ \frac{1}{[L:K]}\sum_{\sigma}  \sum_{\gb}
\sum_{\ga}
 \frac{{}^{\gb^{-1}}\lambda(\gb) r_\ga {}^{(\ga\gb)^{-1}}\lambda(\ga)}{{}^\sigma \psi'(\gb)}  \,{}^\sigma f} \\[10pt]
& = \displaystyle{ \frac{1}{[L:K]}\sum_{\sigma}    \sum_{\ga}
r_\ga g_{\sigma}(\lambda_{{}^{\ga^{-1}}\lambda(\ga)}) {}^{\ga^{-1}}\lambda(\ga) \,{}^\sigma f} \\[10pt]
& = \displaystyle{ \frac{1}{[L:K]}\sum_{\sigma}    \sum_{\ga}
r_\ga {}^\sigma {\psi'(\ga) }g_\sigma(\lambda) \,{}^\sigma f} \\[10pt]
& = \displaystyle{ \frac{1}{[L:K]} \Psi^* \left(\sum_{\sigma}
g_\sigma(\lambda) \,{}^\sigma f \right) = \Psi^*(h)}\,.\qquad \Box
\end{array}
$$

The transitivity of the action of $\iota(E^*)$ on the set of
elliptic directions follows from the equality~(\ref{igualtat}). To
finish the proof of Theorem~\ref{modular}, it remains to determine
the $q$-expansions of the normalized elliptic directions. For it,
first we need a technical lemma.

\begin{lemma}\label{tecnic}
Let $\ell\colon I(\gm) \to L^*$ be a map such that
$\ell(\ga)=\psi(\ga)$ for all $\ga=\operatorname{id}$ in
$\Gal(\overline{L}/K)$. Let $\tau\colon \Gal(\overline{L}/K) \to
\Gal(L/K)$ be a map such that $\ell(\ga \gb)=\ell(\ga)
{}^{\tau(\ga)}\ell(\gb)$ for all $\ga\in I(\gm)$. Then, the
identity
\begin{equation}\label{tecnic1}
\frac{1}{[L:K]} \sum_{\sigma \in \Phi} \beta_\sigma {}^\sigma
\psi(\gc) = \ell(\gc)
\end{equation}
holds for all $\gc\in I(\gm)$ if and only if
\begin{equation}\label{tecnic2}
\beta_\sigma = \sum_{\ga\in \Gal(\overline{L}/K)}
\frac{\ell(\ga)}{{}^\sigma \psi(\ga)} \quad \text{ and } \quad
\sum_{\sigma\in \Phi} \beta_\sigma = [L:K]\,.
\end{equation}
\end{lemma}

\noindent {\bf Proof.} Assume (\ref{tecnic2}). For every $\gc\in
I(\gm)$, we have
$$
\begin{array}{l@{\,=\,}l}
\displaystyle{\sum_{\sigma\in\Phi}\left(\sum_{\ga\in
\Gal(\overline{L}/K)} \frac{\ell(\ga)}{{}^{\sigma}\psi(\ga)}
\right) {}^{\sigma}{\psi(\gc)}} &
\displaystyle{\sum_{\sigma\in\Phi}\left(\sum_{\ga\in\Gal(\overline{L}/K)
} \frac{\ell(\ga\gc )}{{}^{\sigma}\psi(\ga \gc)} \right)
{}^{\sigma}{\psi(\gc)}} \\[8pt]
&
\displaystyle{\ell(\gc)\sum_{\sigma\in\Phi}\left(\sum_{\ga\in\Gal(\overline{L}/K)}
\frac{{}^{\tau(\gc)}\ell(\ga )}{{}^{\sigma}\psi(\ga)} \right)} \\[8pt]
& \displaystyle{\ell(\gc)
\!\!\!\!{\phantom{\int}}^{\tau(\gc)}\left(
\sum_{\ga\in\Gal(\overline{L}/K)}\ell(\ga )
\left(\sum_{\sigma\in\Phi} \frac{1}{{}^{\sigma}\psi(\ga)}\right)\right) } \\[12pt]
& \displaystyle{\ell(\gc)
\!\!\!\!{\phantom{\int}}^{\tau(\gc)}\left(\sum
_{\sigma\in\Phi}\beta_\sigma\right)=\ell(\gc)\,[L:K]\,.}
\end{array}
$$
Now, suppose (\ref{tecnic1}). Fix $\nu\in\Phi$. Note that for
$\sigma\in\Phi$, the characters  $\chi_\sigma$ and $\chi_\nu$ are
equal if and only if $\sigma=\nu$. For every
$\ga\in\Gal(\overline{L}/K)$, one has:
$$
\frac{\ell(\ga)}{{}^{\nu}\psi(\ga)}=\frac{1}{[L:K]}\left(
\beta_\nu+\sum_{\sigma\in\Phi\backslash\{\nu\}}\beta_\sigma
\frac{{}^{\sigma}\psi(\ga)}{ {}^{\nu} \psi(\ga)}
\right)=\frac{1}{[L:K]}\left(
\beta_\nu+\sum_{\sigma\in\Phi\backslash\{\nu\}}\beta_\sigma
(\chi_\sigma \chi_\nu^{-1})(\ga) \right)\,.
$$
Summing over all $\ga$, then
$$
\sum_{\ga\in\Gal(\overline{L}/K)}\frac{\ell(\ga)}{{}^{\nu}\psi'(\ga)}=\beta_\nu+\frac{1}{[L:K]}
\left(\sum_{\sigma\in\Phi\backslash\{\nu\}}\beta_\sigma\sum_{\ga\in\Gal(\overline{L}/K)}
(\chi_\sigma\chi_\nu^{-1})(\ga) \right)=\beta_\nu\,.
$$
The condition $\sum_{\sigma\in \Phi}\beta_\sigma=[L:K]$ is
obtained by replacing $\ga$ with $\cO$ in
(\ref{tecnic1}).\hfill$\Box$

\begin{proposition}\label{p2}
Assume that $\lambda\in [A]$ satisfies
$\operatorname{tr}_\Phi(\lambda)=[L:K]$. Consider the normalized
cusp form
$$h=\frac{1}{[L:K]}\sum_{\sigma\in \Phi} g_\sigma(\lambda) \cdot {}^\sigma f\,.$$
Then,
\begin{itemize}
\item [(i)] one has $$h =\left\{\begin{array}{ll}
\displaystyle{\sum_{(\ga,\gm)=1}{}^{\ga^{-1}}\lambda(\ga) \,q^{\Norm(\ga)}} &  \text{, if $K\nsubseteq E_f$;}\\[5 pt]
\displaystyle{\sum_{(\ga,
\gm)=1}\frac{\Norm(\ga)}{\lambda(\overline \ga)} \,q^{\Norm(\ga)}}
&  \text{, if $K\subseteq E_f$;}
 \end{array}\right. $$
\item [(ii)] for all $\gc\in I(\gmc)$, we have $\iota(\psi'(\gc))^*(h)={}^{\gc^{-1}}\lambda(\gc){}^{\gc^{-1}}h$.
\end{itemize}
\end{proposition}

\noindent {\bf Proof.} For all $\ga\in I(\gm)$, set
$$
\ell(\ga)=\left \{ \begin{array}{ll}{}^{\ga^{-1}}\lambda(\ga) & \text{, if $K\nsubseteq E_f$;} \\[3pt]
\displaystyle{\frac{\Norm(\ga)}{\lambda(\overline\ga)}}& \text{,
if $K\subseteq E_f$.}
\end{array}\right.
$$
It is clear that $\ell(\ga\gb)$ is
$\ell(\ga){}^{\ga^{-1}}\ell(\gb)$ or $\ell(\ga){}^{\overline
{\ga}}\ell(\gb)$ depending on whether $K\nsubseteq E_f$ or not,
respectively.  Since for the case $K\subseteq E_f$ one has
$$
\frac{\ell(\ga^{-1})}{{}^\sigma\psi(\ga^{-1})}=\frac{{}^{\overline
\ga^{-1}}(1/ \ell(\ga))}{{}^\sigma\psi(\ga^{-1})}=
\frac{{}^{\overline\ga^{-1}}
(\Norm(\ga)/\ell(\ga))}{\Norm(\ga)/{}^\sigma\psi(\ga)} =
\frac{{}^{\overline{\ga}^{-1}}\lambda(\overline \ga)}{{}^\sigma
\psi'(\overline\ga)}\,,$$ for all $\sigma\in \Phi$, then in both
cases it follows that $g_\sigma(\lambda)
=\sum_{\ga\in\Gal(\overline{L}/K)} \ell(\ga)/{}^\sigma\psi(\ga)$.
By using Lemma~\ref{tecnic},  a case-by-case computation shows
that for all $\ga\in I(\gm)$ and $\gc\in I(\gmc)$ one has
\begin{equation}\label{miraquebe}
\frac{1}{[L:K]} \sum_{\sigma\in\Phi} g_\sigma(\lambda) {}^\sigma
\psi(\ga){}^\sigma
\psi'(\gc)={}^{\gc^{-1}}\lambda(\gc){}^{\gc^{-1}}\ell(\ga)\,.
\end{equation}
Plugging $\gc=1$ in (\ref{miraquebe}) it follows part (i). Part
(ii) follows from part (i) and (\ref{miraquebe}). \hfill$\Box$
\vskip 0.3truecm

Now, Theorem \ref{modular} in the Introduction follows from
Propositions \ref{p1}, \ref{p3} and \ref{p2}. Note that due to
Proposition \ref{projector}, all one-cocycles in $[A]$ are modular
if and only if $[E:K]=[L:K]$; i.e., when $A$ is $K$-isogenous to
$\Res_{L/K}(C)$. In general, in order to determine a modular
one-cocycle in~$[A]$ a strategy emerges from the previous results.
Indeed, first one can build a one-cocycle~$\lambda\in [A]$ by
solving and combining norm equations. If
$\operatorname{tr}_\Phi(\lambda)\neq 0$, then
$\lambda_{\operatorname{tr}_\Phi(\lambda)}$ is modular since its
$\Phi$-trace equals $[L:K]$. Alternatively, if
$\operatorname{tr}_\Phi(\lambda)= 0$ or in any circumstance, the
nullspace of the $K$-linear map $\operatorname{pr}-[L:K]
\,\operatorname{Id}$ provides all $u\in L$ such that $\lambda_u$
is modular.

\vskip 0.2truecm

We also remark that for the case $K\subseteq E_f$, there are
elliptic quotients of $A_f$ that do not factor through neither~$A$
nor $\overline{A}$. These quotients can be obtained using the
above results plus the Weil involution acting on $A_f$.

\vskip 0.2truecm

We conclude this section with three open questions: one concerning
about the isomorphism $\iota \colon E \rightarrow \End_K^0(A)$ and
the others about the elliptic optimal quotients of $A$. All the
results of the paper hold when we replace~$J_1(N)$ with~$\Jac
(X_\Gamma)$, where~$\Gamma$ is an intermediate congruence subgroup
between~$\Gamma_1(N)$ and~$\Gamma_0(N)$ such that~$f\in
S_2(\Gamma)$ and~$X_\Gamma$ is the modular curve attached to this
subgroup. Although the optimal quotient~$A$ of~$A_f^{(\Gamma)}$
does depend on~$\Gamma$, it is known that~$\iota (T_p)\in
\End_\Q(A_f^{(\Gamma)} )$ and, thus, $\iota(T_p)$ belongs to~$\End
_K(A)$ for all~$\Gamma$.

\begin{question}
Is $\iota (\psi(\ga))\in \End_K(A)$ for all integral ideals $\ga$
and all $\Gamma$?
\end{question}

We ask ourselves whether the $j$-invariants of optimal modular
parametrizations of CM elliptic curves are not far from being also
{\sl optimal} in the sense of having CM by the maximal order
of~$K$. Of course, if $\iota(\cO_K)\subset \End_K(A)$ all optimal
elliptic quotients have multiplication by~$\cO_K$. If $\iota (\eta
(\ga))\in \End_K(A)$ for all integral ideals $\ga\in I(\gmc)$,
then the $j$-invariants of all optimal elliptic quotients are in
the Hilbert class field~$H$. From Cremona's tables ($N<130000$),
we have checked that all optimal elliptic quotients over~$\Q$
with~CM of~$J_0(N)$ have complex multiplication by~$\cO_K$. Also,
the same experimental result has been obtained in all examples
over $\Qbar$ collected by the authors.

\begin{question}
Assume that $\pi\in \Hom_L(A,C)$ is optimal. Does   $C$ have complex
multiplication by $\cO_K$?
\end{question}

\vskip 0.2truecm

And the last question is related to the above Remark~\ref{global}.

\begin{question}
Is it true that the existence of an optimal elliptic quotient of
$A$ having global minimal model over $L$ is equivalent to the
existence of a modular one-cocycle $\lambda\in [A]$ with values
$\lambda(\ga)$ in the ring of integers $\cO_L$ for all integral
ideals $\ga\in I(\gmc)$?
\end{question}

In the next sections, we apply the above results and focus our
attention on Gross's elliptic curves $A(p)$. We also give a
positive answer to the second question mentioned above for the
particular case of level $N=p^2$.

\section{CM elliptic optimal quotients of $J_1(p^2)$}
In the sequel $p$ is a prime $>3$ and such that  $p\equiv 3 \bmod
4$. The discriminant of~$K=\Q(\sqrt{-p})$ is~$-p$. Set $\gp
=\sqrt{-p}\,\cO_K$. Let ${\mathcal X}$ denote the set of Hecke
characters mod~$\gp$ and let ${\mathcal Y}$  be
 the set of Dirichlet characters
$ \eta: {(\cO_K/\gp)}^* \to \C^*$  \text{ such that } $
\eta(-1)=-1$.

To every Hecke character $\psi\in {\mathcal X}$, we attach its
eta-character $\eta$ in ${\mathcal Y}$ defined as in Section 3 by
$ \eta(a) = {\psi((a)) }/{a}$, and it can be easily seen that this
map ${\mathcal X}\to {\mathcal Y}$ is surjective. The  Nebentypus
$\varepsilon: {(\Z/p\Z)}^*\to \C^*$ of the newform $f\in
S_2(\Gamma_1(p^2))$ associated with~$\psi$ is given by
$\varepsilon(n)=\chi(n)\eta(n)$, where~$\chi$ is the quadratic
Dirichlet character associated with~$K$. In this case, we have
that $\ord \varepsilon= (\ord\eta)/2$.

By the results in Section~3, we know that the elliptic optimal
quotients of the abelian variety $A_f$ are defined over a number
field $L$,  which is a cyclic extension of~$H$ of degree $\ord
\varepsilon$   contained in $K_\gp$.

\begin{proposition}\label{rayclassfield} The ray class field $K_\gp$
satisfies~$[K_\gp: H]=(p-1)/2$ and we have~$K_\gp=H\cdot \Q
(\zeta_p)$, where $\zeta_p=e^{2\pi i/p}$.
\end{proposition}

\noindent {\bf Proof.} From the exact sequence
$$ 1\longrightarrow (\cO_K/\gp)^{*}/{\cO_K}^{*}\longrightarrow I(\gp)/P_{1}(\gp)
\longrightarrow I(\cO_K)/P(\cO_K)\longrightarrow 1\,, $$ we know
that the Galois group $\Gal (K_\gp/H)$ is isomorphic to
$(\cO_K/\gp)^{*}/{\cO_K}^{*}$ and, thus, one has $ [K_\gp: H]=
(p-1)/2$. Consider the morphism $\Phi_\gp: I(\gp) \rightarrow
\Gal(H\cdot \Q (\zeta_p)/K)$ given by the Artin symbol. We claim
that $\Phi_p$ has kernel~$P_1(\gp)$, which implies that $K_\gp
\subseteq H\cdot \Q (\zeta_p)$. Indeed, for any ideal $\ga\in
I(\gp)$, we have that $\Phi_\gp(\ga)$ acts trivially on $H$ if and
only if $\ga\in P(\gp)$, that is $\ga=a\cO$. Moreover,
$\Phi_p(a\cO)$ acts trivially on $\Q(\zeta_p)$ if and only if the
Artin symbol $\left (\frac{\Q(\zeta_p)/\Q}{\Norm (a)}\right)$ is
the identity; i.e., $\Norm (a)\equiv 1\pmod {\gp}$ which is
equivalent to $\ga\in P_1(\gp)$ since $\Norm (a)\equiv a^2 \pmod
{\gp}$. Finally, for any subfield $F$ of $\Q(\zeta_p)$ which
contains $K$ we have that $H\cap F=K$ since either $F=K$ or $F/K$
is ramified at~$\gp$. Hence, one has the equality $[H\cdot \Q
(\zeta_p): H]=(p-1)/2=[K_\gp: H]$ and the statement follows.
\hfill $\Box$

\bigskip

We shall need the following lemma.

\begin{lemma}\label{sumzero}
Let $\psi\in \mathcal X$ and denote by $\eta$ and $f$ its
eta-character and newform, respectively. Then,
\begin{itemize}
\item[(i)] for every ideal $\ga\in I(\gp)$, one has
$$
\operatorname{Tr}_{E/K} \left( \psi(\ga) \right) =
\begin{cases}
\displaystyle{ a \sum_{\sigma\in \Phi} {}^\sigma \eta(a)} &
\text{, if $\ga=a\cO_K$;} \\[6pt]
0 & \text{, if $\ga\not\in P(\gp)$;} \\
\end{cases}
$$
\item[(ii)] let $\eta'$ and $f'$ denote the
eta-character and newform associated with $\psi'\in \mathcal X$.
Then, $f'={}^\sigma f$ for some $\sigma \in \Gal(\Qbar/K)$ if and
only if $\ker\eta'=\ker\eta$.
\end{itemize}
\end{lemma}

\noindent {\bf Proof.} First, let us prove (i). When $\ga=a\cO_K$,
the claim on the trace is clear since ${}^\sigma\psi((a))=a
{}^\sigma\eta(a)$. Suppose that $\ga\not\in P(\gp)$, and let $n$
be the order of~$\ga$ in $I(\gp)/P(\eta)$. Notice that $n>1$ and
$\psi(\ga)\not\in K$. For every $\sigma\in \Phi$, we have
${}^\sigma\psi(\ga)=\psi(\ga)\zeta_\sigma$ for some
$\zeta_\sigma\in \mu_n$, where $\mu_n$ denotes the group of
$n\,$th roots of unity. Thus, we have
$$
\sum_{\sigma\in \Phi} {}^\sigma \psi(\ga) = \psi(\ga)
\sum_{\sigma\in \Phi} \zeta_\sigma \in K \,.
$$
Therefore, either $\operatorname{Tr}_{E/K} \left( \psi(\ga)
\right)=0$ or $\psi(\ga)\in K(\mu_n)$. Let us see that the last
possibility does not occur. For it, assume that $\psi(\ga)\in
K(\mu_n)$ which implies that the extension $K(\psi(\ga))/K$ is
normal. Since $n$ is the minimum positive integer such that
$\psi(\ga)^n\in K$, it follows that either $\mu_n\subset K$ or
$\psi(\ga)^{2n}\in K^n$ (see Proposition~2 in \cite{Roh}). Since
$\psi(\ga)\not \in K$, we must have that
$\psi(\ga^{2n})=b^n=\psi((b\cO_K)^{n})$ for some $b\in K$ and,
hence, $\ga^2=b\cO_K$. The class number of $K$ being odd, we get a
contradiction.

Let us prove (ii). If $f'={}^\sigma f$ for some $\sigma \in
\Gal(\Qbar/K)$ then the statement is clear since $\eta'={}^\sigma
\eta$. Now, suppose that $\ker\eta'=\ker\eta$. We claim that
$$
\{ {}^\sigma f \colon \sigma \in \Phi \} \cap \{ {}^\sigma
f'\colon \sigma \in \Phi' \} \neq \varnothing\,,
$$
where $\Phi'$ is the corresponding set of $K$-embeddings
$\Q(\psi')\hookrightarrow \C$. Let us consider the normalized cusp
forms

$$
\begin{array}{l@{\,=\,}l}
h & \displaystyle{\frac{1}{|\Phi|} \sum_{\sigma\in\Phi} {}^\sigma f = q +\dots} \,, \\[18pt]
h' & \displaystyle{\frac{1}{|\Phi'|} \sum_{\sigma\in\Phi'}
{}^\sigma f' = q +\dots}
\end{array}
$$
in $S_2(\Gamma_1(p^2))^{\operatorname{new}}$. Since $K\nsubseteq
\Q(\im \eta)$ and $\ker \eta'=\ker \eta$, there is $\tau\in \Phi$
such that ${}^\tau\eta(a)=\eta'(a)$ for all $a\in \cO_K$ coprime
with~$\gp$. By applying (i), we obtain the equality
$$
h = \sum_{\ga \in P(\gp)} \frac{\operatorname{Tr}_{E/K}
(\psi(\ga))}{|\Phi|} \, q^{\Norm(\ga)} = \sum_{\ga \in P(\gp)}
\frac{\operatorname{Tr}_{E/K} (\psi'(\ga))}{|\Phi'|} \,
q^{\Norm(\ga)} =h'\,.
$$
Therefore, the $\Qbar$-vector spaces generated by $\{{}^\sigma
f\colon \sigma\in \Phi \}$ and $\{{}^\sigma f'\colon \sigma\in
\Phi' \}$ have a common non-zero cusp form, which implies that
$f'={}^\sigma f$ for some $\sigma\in \Gal(\Qbar/\Q)$ (cf.
Proposition~3.2 in \cite{bggp}). Since $h\in \langle {}^\tau
f\colon \tau\in \Gal(\Qbar/K) \rangle \cap \langle {}^\tau
f'\colon \tau\in \Gal(\Qbar/K) \rangle$, it follows that $\sigma
\in\Gal(\Qbar/K)$. \hfill$\Box$

\begin{proposition}\label{dim}
For every positive divisor $d$ of $(p-1)/2$ there is a unique
abelian variety $A_f$ of CM elliptic type of level $p^2$ such that
the Nebentypus of $f$ has order $d$; this abelian variety
satisfies that $K\not\subseteq E_f$ and   $\dim A_f=[H:
K]\varphi(d)$, where $\varphi$ is the Euler function.

\end{proposition}
\noindent{\bf Proof.} Let $d$ be a divisor of $(p-1)/2$ and take
$\psi\in {\mathcal X}$ such that its  eta-character has  order
$2\,d$. Let us denote by $f$ the newform attached to $\psi$, whose
Nebentypus $\varepsilon$ has order $d$. First, let us show that
$K\not\subseteq E_f$. Indeed, let $\psi_c\in \mathcal X$ defined
by $\psi_c(\ga)=\overline {\psi(\overline \ga)}$. The
eta-character and the normalized newform attached to $\psi_c$ are
clearly $\overline \eta$ and $\overline f$, respectively. Since
$\ker \overline\eta=\ker\eta$, Lemma~\ref{sumzero}~(ii) ensures
that $\overline f\in \{ {}^\sigma f\colon \sigma \in \Phi \}$,
which implies $K\not\subseteq E_f$. The same argument can be
applied to another newform $f'$ obtained from $\psi' \in {\mathcal
X} $ whose associated  character $\eta'$ has order $2d$ to show
that $f'$ belongs to $\{ {}^\sigma f\colon \sigma \in \Phi \} $,
which proves that~$A_f$ is unique when  the order of~$\varepsilon$
has been fixed.

Since $K\not\subseteq E_f$, the equality $\dim A_f =[E_f:
\Q]=[E:K]$ holds. Now, we have that $[E\colon K]=| \{{}^\sigma
f\colon \sigma \in \Phi\}| = | \{{}^\sigma \psi\colon \sigma \in
\Phi\} |$. Again using part~(ii) of Lemma~\ref{sumzero}, we obtain
$$[E: K]=
|\{ \sigma \in \Phi\colon \eta={}^\sigma \eta\}|\cdot |\{ {}^\sigma
\eta\colon \sigma \in \Phi\}|=|\{ \sigma \in \Phi\colon
\eta={}^\sigma \eta\}| \cdot \varphi(d)\,.
$$
Since the condition ${}^\sigma \eta=\eta$ is equivalent to
$\psi/{}^\sigma \psi$ being a character of $\Gal (H/K)$,  it follows
$\dim A_f=[H: K]\varphi(d)$. \hfill $\Box$

\begin{remark}
Note that the number of abelian varieties $A_f$ of CM elliptic
type of level~$p^2$ is the number of divisors of $(p-1)/2$. Also
for every number field $L$ intermediate between $H$ and $H\cdot
\Q(\zeta_p)$ there is a unique abelian variety $A_f$ of CM
elliptic type and level~$p^2$ for which $L$ is its splitting field
as defined in Section~3.
\end{remark}

\bigskip

Next, in order  to show that the CM elliptic optimal quotients
of~$A_f$ in $J_1(p^2)$ have endomorphism ring isomorphic to~$\cO_K$,
 we shall need to use some auxiliary congruence
subgroups of~$\SL_2(\Z)$ of level~$p^2$. To this end, fix $f$ a
newform in $S_2(\Gamma_1(p^2))$ attached to a Hecke character
$\psi\in \mathcal X$. Let~$\varepsilon$ denote the Nebentypus of
$f$. Let us consider  the following congruence subgroups of
level~$p^2$:
$$
\begin{array}{lll}
\Gamma_p&=&\left\{ \left( \begin{array}{cc} a & b\\ c & d
\end{array}\right)\in \Gamma_0( p^2)\colon a\equiv d\equiv 1 \pmod p
\right\}\,,
\end{array}$$
and $\Gamma_\varepsilon$ as in the introduction; i.e.,
$$
\begin{array}{lll}
\Gamma_\varepsilon&=&\left\{ \left( \begin{array}{cc} a & b\\ c & d
\end{array}\right)\in \Gamma_0( p^2)\colon \varepsilon (d)=1
\right\}\,.
\end{array}$$
 It is clear that $\Gamma_1(p^2)\subseteq \Gamma_p\subseteq
\Gamma_\varepsilon$ and $f\in S_2(\Gamma_\varepsilon)$. For any
intermediate  congruence subgroup $\Gamma$ of level $p^2$
satisfying $ \Gamma_1(p^2) \subseteq \Gamma \subseteq
\Gamma_\varepsilon\,, $ let $X_{\Gamma}$ be the modular curve
over~$\Q$ attached to $\Gamma$.  We shall denote by
$A_f^{(\Gamma)}$  the optimal quotient of the jacobian
of~$X_{\Gamma}$  attached to~$f$ by Shimura. More precisely,
let~$I_f$ be the annihilator of~$f$ in the Hecke algebra acting
on~$\operatorname{Jac}(X_{\Gamma})$. Then,
$$A_f^{(\Gamma)}=\operatorname{Jac}(X_{\Gamma})/I_f\left(\operatorname{Jac}(X_{\Gamma})\right)\,
.$$

\begin{proposition}\label{optimal}
Let $f$ and $\Gamma$ be as above. Then, all elliptic optimal
quotients of $A_f^{(\Gamma)}$ have complex multiplication
by~$\cO_K$.
\end{proposition}

\noindent{\bf Proof.} Fix an elliptic direction in $\Omega^1(A_f)$
and let $C_\Gamma$ be an elliptic optimal quotient attached to
this direction. By Proposition~\ref{dim} and
Theorem~\ref{sistema}, we know that $K\not\subseteq E_f$ and thus
all endomorphisms of~$A_f^{(\Gamma)}$ are defined over its
splitting field, say $L$, that satisfies $L\subseteq K_\gp$. Let
$c_\Gamma$ denote the conductor of the order
$\cO_\Gamma\simeq\End_L(C_\Gamma)$ in $\cO_K$. We want to show
that $c_\Gamma=1$, and split the proof in three steps.

\medskip

{\it Step 1: $c_{\Gamma}\mid 2$ for all $\Gamma$.} Since $\End
(C_\Gamma)=\End_L(C_\Gamma)$, one has that $L$ contains the ring
class field of $\cO_\Gamma$, say $K_\Gamma$. Notice that $K_\Gamma
\subseteq L\subseteq K_\gp$. But $p\nmid c_\Gamma$, since otherwise
$p$ must divide $[L: H]$ (cf. Proposition~7.24 in~\cite{cox89}) and
this degree is a divisor of $(p-1)/2$. Hence, $K_\Gamma$ is an
unramified extension of the Hilbert class field and, therefore, it
must coincide with~$H$. Again by Proposition 7.24 in \cite{cox89},
we obtain that $c_\Gamma\mid 2$.

\bigskip

{\it  Step 2: $c_\Gamma$ does not depend on  $\Gamma$.} We consider
the natural projection $\pi: X_{\Gamma} \rightarrow
X_{\Gamma_\varepsilon}$. The degree of~$\pi$ is odd since it divides
$ [\Gamma_1(p^2): \Gamma_0(p^2)/\{\pm 1\}]=p( p-1)/2$ and $p\equiv 3
\bmod 4$.

Let $\pi_{\Gamma,\Gamma_\varepsilon}: \Jac (X_\Gamma)\longrightarrow
A_{\Gamma,\Gamma_\varepsilon}$ be the optimal quotient over $\Q$ for
which there is an isogeny $\nu:
A_{\Gamma,\Gamma_\varepsilon}\longrightarrow  \Jac
(X_{\Gamma_\varepsilon})$ defined over $\Q$ rending the following
diagram
$$ \xymatrix{
\Jac (X_\Gamma)\ar[rr]^-{\pi_*}
\ar[rd]_-{\pi_{\Gamma,\Gamma_\varepsilon}} & & \Jac
(X_{\Gamma_\varepsilon}) \\
          & A_{\Gamma,\Gamma_\varepsilon} \ar[ru]_{\nu} &}
$$
commutative. Since every element of the group
$H_1(X_{\Gamma_\varepsilon},\Z)/\pi_*( H_1(X_\Gamma,\Z))$ has order
dividing~$\deg \pi$, the cardinality of this group is odd. From the
group isomorphism $ \ker \nu\simeq
H_1(X_{\Gamma_\varepsilon},\Z)/\pi_*( H_1(X_\Gamma,\Z)) $,  it
follows that  $\deg \nu $ is odd. Since $A_f^{(\Gamma)}$ is an
optimal quotient of $A_{\Gamma,\Gamma_\varepsilon}$, there is an
isogeny $\nu_f: A_f^{(\Gamma)}\longrightarrow
A_f^{(\Gamma_\varepsilon)}$ whose degree divides $\deg \nu$. Hence,
for every optimal elliptic quotient $\pi_\Gamma:
A_f^{(\Gamma)}\longrightarrow C_{\Gamma}$ there is an optimal
elliptic  quotient $\pi_\varepsilon:
A_f^{(\Gamma_\varepsilon)}\longrightarrow C_{\Gamma_\varepsilon}$
 and an isogeny $\mu: C_{\Gamma}\longrightarrow C_{\Gamma_\varepsilon}$ rending
the diagram
$$ \xymatrix{
A_f^{(\Gamma)} \ar[rr]^-{\nu_f} \ar[dd]_{\pi_\Gamma}  &&
A_f^{(\Gamma_\varepsilon)} \ar[dd]^{\pi_\varepsilon}
\\  \\
C_{\Gamma} \ar[rr]^-{\mu} && C_{\Gamma_\varepsilon}}
$$
commutative. It is clear that $\deg \mu$ is odd since it divides
$\deg \nu_f$. Thus~$c_{\Gamma_\varepsilon}$ and~$c_\Gamma$ can only
differ by an odd factor, which implies that $c_\Gamma$ is
independent of the group  $\Gamma$.
\bigskip

{\it Step 3: $c_{\Gamma}=1$ for all $\Gamma$.} Now, it suffices to
prove $c_\Gamma=1$ for a particular subgroup~$\Gamma$. We consider
$\Gamma=\Gamma_p$. Following Shimura in~\cite{shimuraoptimal}, we
know that  the matrix
$$
\left (\begin{array}{cc} 1& 1/p \\0 & 1
\end{array}\right )
$$
lies in the normalizer of $\Gamma_p$ in $\operatorname{SL}_2(\R)$
and provides an automorphism~$u$ of~$X_{\Gamma_p}$ of order~$p$.
Set
$$
G= \sum_{\stackrel{1\leq i<p}{\chi (i)=1}}(u^*)^i \in \End \Jac
(X_{\Gamma_p})\,.
$$

We claim that $G$ leaves stable the subvariety  $I_f (\Jac
(X_{\Gamma_p}))$, which  is equivalent to saying that $G$ leaves
stable the vector space generated by the set of eigenforms
in~$S_2(\Gamma_p)$ which are not Galois conjugates of $f$. In
fact, the action of~$G$ on all eigenforms of~$S_2(\Gamma_p)$ can
be described as follows. It is well-known that if we denote
by~$\operatorname{New}_\Gamma$ the set of normalized newforms
in~$S_2(\Gamma)$, then the set of normalized eigenforms
in~$S_2(\Gamma_p)$ is the disjoint union
of~$\operatorname{New}_{\Gamma_p}$,  ${\mathcal S_1}$, and
${\mathcal S_2}$, where ${\mathcal
S_1}={\operatorname{New}}_{\Gamma_1(p)} \cap S_2(\Gamma_p)$,
${\mathcal S_2}= B_p({\operatorname{New}}_{\Gamma_1(p)}) \cap
S_2(\Gamma_p)$, and~$B_p$ is the operator acting
as~$B_p(h(q))=h(q^p)$. With~$\zeta_p= e^{2 \pi i/p}$ and from the
equality
$$ \sum_{\stackrel{1\leq i < p}{\chi (i)=1}} \zeta_p^i=\frac{-1+\sqrt{-p}}{2}\,,$$ it can be easily
checked that every eigenform $h(q)=\sum_{n\geq 1}b_n q^n\in
S_2(\Gamma_p)$ satisfies:
$$G^*(h)=\left\{
\begin{array}{ll}
\displaystyle{\frac{-1+\sqrt{-p}}{2}h+ \frac{p-\sqrt{-p}}{2}\,
b_p\, B_p(h)} & \text{, if
$h\in \operatorname{New}_{\Gamma_p}\cup \mathcal S_1$,}\\[10 pt]
\displaystyle{\frac{p-1 }{2}\,h} & \text{, if $h\in {\mathcal S_2}
$.}
\end{array}
\right.
$$
The claim follows from the fact that all $h\in
\operatorname{New}_{\Gamma_p}$  have level $p^2$ and Nebentypus
whose conductor divides $p$ and, thus,  $b_p=0$ (see
subsection~1.8 in~\cite{deligne:serre74}).

Since $G$ leaves stable the subvariety  $I_f (\Jac (X_{\Gamma_p}))$,
then $G$ induces an endomorphism of $A_f^{(\Gamma_p)}$, which we
still denote by $G$. Due to the fact that  $G$ acts on
$\Omega^1(A_f^{(\Gamma_p)})$ as the multiplication by
$(-1+\sqrt{-p})/2$, it follows that $G$ leaves stable all
subvarieties of $A^{(\Gamma_p)}$. Thus, $(-1+\sqrt{-p})/2 \in
\cO_{\Gamma_p}$ and the statement follows.
 \hfill $\Box$

\bigskip

As for Gross's elliptic curves, we obtain the following result,
which concludes the proof of Theorem \ref{Thm4}.

\begin{cor}
Let $f$ be a CM normalized newform with trivial Nebentypus. The
elliptic curve $A(p)$ and its Galois conjugates are the optimal
quotients of $A_f^{(\Gamma)}$ over the Hilbert class field $H$,
for all subgroups~$\Gamma$ with $\Gamma_1(p^2)\subseteq \Gamma
\subseteq \Gamma_0(p^2)$.
\end{cor}

\noindent{\bf Proof.} By Theorem~20.1 in~\cite{gross80}, we know
that $A(p)$ is a quotient of $J_0(p^2)$ defined over~$H$, attached
to a newform $f$ with trivial Nebentypus. Notice that the
corresponding field $L$ coincides with the Hilbert class
field~$H$. Since we have~$K\not\subseteq E_f$, by
Theorems~\ref{teoremaA} and~\ref{sistema}, every elliptic optimal
quotient~$C_\Gamma$ of~$A_f^{(\Gamma)}$ is defined over~$H$ and
the abelian variety~$A_f^{(\Gamma)}$ is simple over~$K$.
Since~$\dim A_f^{(\Gamma)}=[H: K]$, it follows
that~$A_f^{(\Gamma)}$ is $K$-isogenous to the Weil
restriction~$\Res_{H/K} C_\Gamma$. In~\cite{gross80}, Gross shows
that~$A_f^{(\Gamma)}$ is $K$-isogenous to $\Res_{H/K} A(p)$.
Therefore, on the one hand, there is $\sigma\in \Gal(H/K)$ such
that~$A(p)$ and ${}^\sigma C_\Gamma$ are $\Qbar$-isomorphic. On
the other hand, by Theorem~\ref{optimal}, $A(p)$ and ${}^\sigma
C_\Gamma$ are $H$-isogenous. Hence, $A(p)$ is $H$-isomorphic to
${}^\sigma C_\Gamma$ and the claim follows. \hfill $\Box$

\section{Canonical CM elliptic direction for $A(p)$}

When the class number of $K$ is greater than one, there are
infinitely many elliptic directions in~$S_2(\Gamma_0(p^2))$
attached to different parametrizations~$J_0(p^2)\to A(p)$. Here,
we shall emphasize one of them (we call it canonical) in terms of
a particular one-cocycle that can be constructed by means of the
Dedekind eta-function.

Let $\cO_H$ be the ring of integers of the Hilbert class
field~$H$. For all $a\in K$ coprime with $\gp$, we denote by
$(\frac{a}{\gp})$ the Jacobi symbol $(\frac{m}{p})$, where $m$ is
an integer such that $a\equiv m\pmod{\gp}$. One has
$\eta(a)=(\frac{a}{\gp})$. By \cite{gross82}, we know that there
is a unique map $\delta: I(\gp)\rightarrow H$ with the following
two requirements:
\begin{itemize}

\item[(i)] $\delta(\ga)^{12}=\Delta(\cO)/\Delta(\ga)$,
\item[(ii)] $\displaystyle{\left(\frac{\Norm_{H/K}
(\delta(\ga))}{\gp}\right)=1}$,
\end{itemize}
for all $\ga\in I(\gp)$. Moreover, this map also satisfies the
following conditions:

\begin{itemize}
\item [(iii)] $\delta(\ga) \cO_H=\ga\cO_H$,

\item[(iv)] $\delta(\ga\cdot\gb)=\delta(\ga)\cdot{} ^{\ga^{-1}}\delta(\gb)$ for all
$\ga,\gb\in I(\gp)$,

\item[(v)] $\delta(\overline \ga)=\overline{\delta(\ga)}$ for all
$\ga\in I(\gp)$.
\end{itemize}
By taking into account conditions (ii) and (iv), and since $[H:K]$
is odd, we also obtain:

\begin{itemize}
\item[(vi)] for all $\ga\in P(\gp)$,  one has
$\delta(\overline \ga)\in K$ and
$\displaystyle{\left(\frac{\delta(\ga)}{\gp}\right)=1}$.
\end{itemize}

For every $\ga\in I(\gp)$, we set
\begin{equation}\label{lambda2}
\lambda(\ga):={}^\ga \delta(\ga)=\frac{\Norm
(\ga)}{\delta(\overline{\ga})}\,.\end{equation} The map $\lambda:
I(\gp)\longrightarrow H$ also satisfies conditions (ii), (iii), (v),
and (vi). But now conditions (i) and (iv) are replaced with
\begin{itemize}
\item [(i')] $\lambda(\ga)^{12}=\Norm(\ga)^{12} \displaystyle{\frac{\Delta({\overline{\ga}})}{\Delta(\cO_K)}}$,
\end{itemize}
and the one-cocycle condition:
\begin{itemize}
\item [(iv')] $\lambda(\ga\cdot\gb
)=\lambda(\ga)\cdot{}^\ga\lambda(\gb)$, for all $\ga$, $\gb \in
I(\gp)$.
\end{itemize}
Conditions (vi) and (iv') imply that the one-cocycle $\lambda$
belongs to $[A_f]$ for all $A_f$ of CM elliptic type and
level~$p^2$.
\begin{remark}
Notice that the above one-cycle $\lambda$ can be effectively
computed by using the Dedekind eta-function on ideals (as
Rodríguez-Villegas does in~\cite{Villegas91}), and it concides
with what Hajir denotes~$\phi$ in~Definition~2.3 in~\cite{hajir}.
\end{remark}

Let $f$ denote the normalized newform in $S_2(\Gamma_0(p^2))$
attached to a Hecke character~$\psi$ whose eta-character has order
$2$. By Section~3, the splitting field~$L$ of~$A_f$ is~$H$. Let
$S_2(A_f)$ be the $\C$-vector space generated by the Galois
conjugates of the newform~$f$ attached to $\psi$ and let $\omega$
denote a N\'eron differential of Gross's elliptic curve $A(p)$.

\begin{proposition}\label{direction} Let $f$ be as above.  There is an optimal quotient
$\pi: J_0(p^2)\longrightarrow A(p)$ such that $\pi^*(\omega)= c\,
g(q)\,dq/q$, where
$$g(q)= \sum_{(\ga,\gp)=1} \delta(\ga) q^{\Norm (\ga)}\in S_2(A_f)\,,$$
and $c\in \Z$ is a unit in $\Z[\frac{1}{2p}]$.

\end{proposition}

\noindent{\bf Proof.} By Lemma \ref{dim},  we have  $[E:K]=[L:K]$
and, thus,  all one-cocycles in $[A_f]$ are modular. Therefore, by
Theorem \ref{modular} we have that
$$g(q)=\sum_{(\ga,\gp)=1} {}^{\ga^{-1}}\lambda(\ga) q^{\Norm (\ga)}=\sum_{(\ga,\gp)=1} \delta(\ga) q^{\Norm (\ga)}$$
is a normalized cusp form in $S_2(A_f)$ for which there is an
optimal elliptic quotient $ C_\lambda$ given by  the lattice
$$\Lambda_g=\left\{ 2\pi\,i \int_\gamma g (z) dz\colon \gamma\in
H_1(X_0(p^2),\Z)\right\}\,.$$ Since $\delta(\overline
\ga)=\overline{\delta(\ga)}$ for all $\ga\in I(\gp)$, it follows
that $g(q)\in H_0[[q]]$. Thus, $g (q)\,dq/q\in
\Omega^1(X_0(p^2))_{/H_0}$. Hence, the natural morphism $\pi:
X_0(p^2)\rightarrow C_\lambda$ is defined over~$H_0$. Notice that
necessarily one has $\Lambda_g=\Omega\cdot \cO_K$, for some
$\Omega\in\C^*$. Indeed, $\cO_K$ is the only ideal $\ga$ such that
$j(\ga)=\overline{j(\ga)}$ since~$[H:K]$ is odd. Thus, we have
$j(\Lambda_g)=j(\cO_K)$. Since~$A_f$ is $\Q$-isogenous to
$\Res_{H_0/\Q}(C_\lambda)$ and to~$\Res_{H_0/\Q}(A(p))$, it
follows that~$C_\lambda$ and~$A(p)$ are~$H_0$-isogenous and,
therefore, $H_0$-isomorphic. Therefore, there exists $c\in H_0^*$
such that $\pi^* (\omega)=c\, g(q)\,dq/q$. It is clear that
$\Delta(\Lambda_g)=-p^3 c^{12}$.

The Manin ideal attached to  $\pi$ is $ c\, \cO_{H_0}$ (we refer
to Section~4 in~\cite{gola01} for more details on the Manin
ideal). By Propositions~4.1 and~4.2 in~\cite{gola01}, we know that
$c \cO_{H_0}$ is an integral ideal and it can only be divided by
primes lying over~$2$ or~$p$. Now, we want to prove that $c\in\Z$.
Since $\pi^*(\omega/c)=g \,dq/q$, the one-cocycle attached
to~$\omega/c$ is~$\lambda$. This means that for every $\ga\in
I(\gp)$ there is an isogeny of degree~$\Norm(\ga)$,
$${\mu}\colon {}^{\ga^{-1}} C_\lambda\rightarrow
C_\lambda\,,$$ such that $\mu ^*(\omega/c)=
{}^{\ga^{-1}}\lambda(\ga)\cdot {}^{\ga^{-1}}(\omega/c)$. Taking into
account that $j(\ga)={}^{\ga^{-1}} j(\cO_K)$, we obtain that the
lattice corresponding to ${}^{\ga^{-1}} C_\lambda$ is
$\frac{1}{\delta(\ga)} \cdot \Omega\ga $. Finally, we have that:
$${}^{\ga^{-1}} \Delta (\Omega\cO_K)=\Delta\left(\frac{1}{\delta
(\ga)} \Omega \ga\right)={\delta(\ga)^{12}} \Delta(\Omega\ga)=
\frac{\Delta (\cO_K)}{\Delta({\ga})}{\Delta(\Omega
\ga)}=\Delta(\Omega\cO_K)\,.
$$
Therefore, $\Delta(\Lambda)\in K\cap H_0=\Q$ and $c^{12}\in\Q$.
Since $\Q(c)\subseteq H$ is unramified outside~$p$ and there is
not a real quadratic field of discriminant $p$, it follows that
$c^3\in \Q$. Finally, since~$H$ does not contain the $3$rd roots
of unity (recall $p>3$), one obtains $c\in \Q$. \hfill $\Box$

\begin{remark}
Since $c\in K^*$, the one-cocycle attached to $\omega$ is also
$\lambda$.  In this sense, we say that the normalized cusp form
$g$ is the canonical cusp form attached to $A(p)$.
\end{remark}

For when the class number of $K$ is $1$ (that is, $p=7$, $ 11$,
$19$, $43$, $67$, $163$), one has that $\pi$ is defined over~$\Q$
and $c$ coincides with the (classical) Manin constant. Then, $c=\pm
1$ in these cases since Manin's conjecture has been checked for all
elliptic curves over~$\Q$ with conductor $\leq 130000$ in Cremona's
tables. We have computed $c$ for the remaining primes $p\leq 100$
(that is, $p= 23$, $31$, $47$, $59$, $71$, $79$, $83$) and we have
also obtained that $c=\pm 1$. It seems reasonable to expect $c=\pm
1$ for all~$A(p)$.

\begin{remark}In general, as already mentioned, there are infinitely many
normalized cusp forms $g'\in S_2(A_f)$ whose directions are
pullbacks of $\Omega^1(A(p))$ under modular para\-me\-tri\-za\-tions
$\pi': A_f\rightarrow A(p)$. For each one of them, there is a
one-cocycle $\lambda'$ (cohomologous to $\lambda$) such that
$$g'=\sum_{\ga} {}^{{\ga}^{-1}}\lambda'(\ga) q^{\Norm (\ga)}\in S_2(A_f)\,,$$
and a constant $c'\in H_0$ with $\pi'^*(\omega)=c'g'$. The concern
on whether the constant $c$ is $\pm 1$ is already in
\cite{gross80}, see Question~23.2.2 in pag.\,81, but without
fixing~$\pi'$. However, $c'\neq \pm 1$ unless $\pi'=\pi$
(canonical) as in Theorem~\ref{direction}, although the Manin
ideal attached to any~$\pi'\neq \pi$ might still be~$\cO_K$ as
well.
\end{remark}

We end this section giving an expression for  the transcendental
$\Omega\in \C^*$ attached to the lattice $\Lambda$ of $A(p)$,
which generalizes the one given by Gross in \cite{gross80} for
when $K$ has class number one. Keeping the above notations, as
in~\cite{hajir}, we set
$$
\rho:=\prod_{ \stackrel{\gb\in\Gal (H/K)}{(\gb,\gp)=1}}
\frac{\delta(\gb)}{\psi(\gb)}\,.
$$
It is clear that $\rho$ is well-defined, independent of the Galois
conjugate of $\psi$,  and $\rho\in \cO_H^*$. Let $h$ denote the
class number of~$K$, and consider
$$\{\cO_K,\gb_1,\cdots,\gb_{(h-1)/2},\cdots,\overline\gb_1,\cdots,\overline\gb_{(h-1)/2}\}$$
a set of representatives of $\Gal(H/K)$ with $(\gb_i,\gp)=1$.
Then, we can rewrite
\begin{equation}\label{conj}
\rho=\prod_{i=1}^{(h-1)/2}\frac{ \delta(\gb_i)\,\delta(\overline
\gb_i)}{\Norm (\gb_i)}=\prod_{\gb\in\Gal (H/K)}\frac{
\delta(\gb)}{\sqrt{\Norm (\gb)}}\,.
\end{equation}
Indeed, since $\delta(\gb)/\psi(\gb)$ is independent of the class of
$\gb$ in $\Gal(H/K)$, it suffices to prove that
$\psi(\gb)\cdot\psi(\overline \gb)=\Norm (\gb)$. But this is a
consequence of
$$
\left(\frac{\Norm (\gb)}{p}\right)=\left(\frac{\Norm
(\gb)}{p}\right)^h=\left(\frac{\beta}{\gp}\right)\left(\frac{\overline\beta}{\gp}\right)=\left(\frac{\beta}{\gp}\right)^2=1\,,
$$
where $\beta\in K$ is a generator of $\gb^h$. Observe that $\rho$
is a positive unit in $\cO_{H_0}^*$.

\begin{proposition}\label{trascendent}
Let $\Lambda=\Omega\cdot\cO_K$ be the lattice attached to~$A(p)$.
Then,
$$\Omega=\pm \, i^{(p+1)/4}\,\,
\displaystyle{\sqrt[\leftroot{6}\uproot{32} h]{ \displaystyle{
{\rho\cdot(2\pi)^{(2 h+1-p)/4}\cdot\sqrt{p}^{(1-3h)/2}}}\cdot
\displaystyle{\prod_{\stackrel{1\leq
m<p}{\chi(m)=1}}\Gamma\left(\frac{m}{p}\right)}}}  \,, $$ where
the $h$-th root is taken to be real.
\end{proposition}

\noindent{\bf Proof.} By the Chowla-Selberg formula~\cite{chowla},
we know that
$$ \prod_{\ga \in\Gal(H/K)}  \Norm (\ga)^{-6}  \Delta
(\tau_\ga)=\left(\frac{2\,\pi}{p}\right)^{6 h}\left(
\prod_{m=1}^{p-1}\Gamma\left(\frac{m}{p}\right)^{\chi(m)}\right)^6\,,
$$
where $\langle 1,\tau_\ga\rangle=\frac{1}{\Norm(\ga)}\ga$. Since
$\lambda$ is the one-cocycle attached to $\omega$, we have that:
\begin{equation}\label{lattice}
\Delta (\tau_\ga)=\Norm (\ga)^{12} \Delta(\ga)=\Norm(\ga)^{12}
\Delta\left( \frac{\Omega}{\delta(\ga)}\,\ga\right)
\frac{\Omega^{12}}{\delta(\ga)^{12}}=-p^3\frac{\Norm
(\ga)^{12}}{\delta(\ga)^{12}} \, \Omega^{12}\,.
\end{equation}
Combining (\ref{conj}), (\ref{lattice}), and Gauss's identity
$$\prod_{i=1}^{n-1}\Gamma\left(
\frac{i}{n}\right)=(2\pi)^{(n-1)/2} n^{-1/2}\,,
$$
the statement follows by taking into account that $\Omega$ lies in
$\R$ or $i\,\R$ according to $p\equiv -1 \pmod 8$ or not (cf.
\cite{gross80}). \hfill $\Box$

\bigskip

As by result, we obtain the following fact, which concludes the
proof of Theorem~\ref{Thm5}.

\begin{cor} With the above notations, one has
$$\left\{ 2\pi\,i \int_\gamma g (z) dz\colon \gamma\in
H_1(X_0(p^2),\Z)\right\} = \frac{1}{c} \cdot \Omega \cdot \cO_K
\,.
$$
\end{cor}

\section{CM elliptic directions for non-trivial Nebentypus}

In this section, we shall consider arbitrary Hecke characters
mod~$\gp$. Let $\psi$ in~$\mathcal X$ and let~$\eta$ be its
eta-character. Let $f$ denote the normalized newform attached
to~$\psi$. In order to find the elliptic directions in $S_2(A_f)$,
one needs to determine the  modular one-cocycles $\lambda_u$ in
$[A_f]$. Then, the normalized cusp forms
$$
g_u = \sum_{(\ga,\gp)=1} {}^{\ga^{-1}}\lambda_u(\ga)\,q^{\Norm(\ga)}
$$
are the elliptic directions in $S_2(A_f)$. Recall that in the
particular case $\eta^2=1$, all one-cocycles are modular. In
general, as explained above, to find the modular one-cocycles
amounts to an eigenvector problem. In our particular setting, the
following lemma will be useful since it will allow to handle
certain linear systems by means of a quotient polynomial ring.

\begin{lemma}\label{general}
Let $M/F$ be a cyclic field extension of degree~$k$. Fix~$\tau$ a
generator of~$\Gal (M/F)$, and let $\mu_k$ be the group of $k\,$th
roots of unity. Let $\cE=\End_{F[\Gal(M/F)]}(M)$ be the
$F$-algebra of $\Gal(M/F)$-equivariant $F$-linear
endo\-mor\-phisms of $M$. One has:
\begin{itemize}
\item[(i)]
the  map $\Theta: F[X]/(X^k-1)\longrightarrow \cE$ given by
$$ \Theta (\sum_{i=1}^k a_i\, X^i)(u)=\sum_{i=1}^k a_i\,{}^{\tau^i}u\,, \quad \text{for all }u\in M\,,$$
is well-defined and an isomorphism of $F$-algebras.

\item[(ii)]
For every $p(X)\in F [X]/(X^k-1)$, let $\cZ=\{ \zeta\in\mu_k\colon
p(\zeta)=0\}$. Then, the endomorphism $G=\Theta(p(X))$
diagonalizes and its characteristic polynomial is
$$(-1)^k\prod_{i=1}^k\left( X-p(\zeta_k^i)\right) \,,$$
where $\zeta_k=e^{2\pi i/k}$. We have $ \dim_{F} \ker G=|\cZ|$,
and
\begin{equation}\label{ker} \ker G=\Theta \left(
\frac{X^k-1}{\displaystyle{\prod_{\zeta\in\cZ}(X-\zeta)}}
\right)(M)\,.\end{equation}
\end{itemize}
\end{lemma}

\noindent{\bf Proof.}   It is obvious that $\Theta$ is
well-defined and a morphism of $F$-algebras.  Chose $\alpha\in M$
such that $\{ {}^{\tau^i}\alpha\}_{1\leq i \leq k}$ is a $F$-basis
of $M$. The morphism $\Theta$ is injective because
$\Theta(q(X))=0$ implies that $\Theta(q(X))(\alpha)=0$ and, then,
$q(X)=0$. For a given $G\in\cE$, we have that $G(\alpha)=\sum
_{i=1}^k a_i {}^{\tau^i}\alpha$ for some $a_i\in F$ and, thus,
$G(u)=\sum _{i=1}^k a_i {}^{\tau^i} u$ for all $u\in M$.
Therefore, $\Theta$ is surjective and  part (i) is proved.

We consider the $F$-algebra monomorphism $\Psi: \cE\longrightarrow
\End_{F} F[X]/(X^k-1)$ defined by $\Psi (G)= \widehat G$, where
\begin{equation}\label{condition}
\widehat{G}(q(X))=\Theta^{-1}(G)\cdot q(X)\,, \quad\text{for all
}q(X)\in F[X]/(X^k-1)\,.
\end{equation}
Now it suffices to prove part (ii)  for $\widehat G$.  Note that
for any field extension~$F_0/F$, the relation (\ref{condition})
allows us to consider $\widehat G$ as a $F_0 $-linear endomorphism
of $F_0[X]/(X^k-1)$.

Let $G= \Theta(p(X))$. The set of eigenvalues of $\widehat G$ is
$\{ p(\zeta_k^i)\colon 1\leq i\leq k\}$. Indeed, if $\beta\in F_0$
is an eigenvalue of eigenvector $q(X)\in F_0[X]/(X^k-1)$, then
there exists $\zeta\in\mu_k$ such that $q(\zeta)\neq 0$ and, thus,
$\beta=p(\zeta)$. Conversely, if $\beta=p(\zeta)$ for some
$\zeta\in\mu_k$ then $q(X)=\prod_{\zeta'\in
\mu_k\backslash\{\zeta\}} (X-\zeta')$ is an eigenvector with
eigenvalue~$\beta$. Notice that all eigenvalues of $\widehat G$
are in $F_0=F(\mu_k)$.

Now, let $\beta=p(\zeta)$ for some $\zeta\in\mu_k$ and we will
prove that
 $$\dim_{F_0} \ker (\widehat{G}-\beta \id)=|\{\zeta\in \mu_k\colon
 p(\zeta)= \beta\}|\,,$$ which implies part (ii) except for the
equality~(\ref{ker}). Note that by a translation of $\widehat G$,
we can (and do) assume $\beta=0$. Then, one has
$$
\begin{array}{l@{\,=\,}l}
\ker \widehat{G} &  \{q(X)\in F_0[X]/(X^k-1)\colon q(\zeta)=0
\text{ for all }\zeta\in\mu_k\backslash \cZ\}= \\[ 8 pt]
 & \{q(X)\in F_0[X]/(X^k-1)\colon q(X)=\displaystyle{\prod_{\zeta\in\mu_k\backslash \cZ}}(X-\zeta)\,r(X)\,, \, \deg
r < |\cZ| \}\,.
\end{array}
$$
It follows that $\dim_{F_0} \ker \widehat G=|\cZ|$ and
 $\ker \widehat G=\ker
(\Psi\circ\Theta)(\prod_{\zeta\in\cZ}(X-\zeta))$. Finally, the
equality (\ref{ker}) is a consequence of the fact that
$q(X)=\prod_{\zeta\in\cZ}(X-\zeta)\in F[X] $ is coprime with
$r(X)=(X^k-1)/p(X)$ and $q(X)\cdot r(X)$ is zero in
$F[X]/(X^k-1)$. \hfill $\Box$

\bigskip

Now, we focus our attention to the Hecke character $\psi\in\mathcal
X$. For the~sake of simplicity, let us assume that its eta-character
satisfies~$\operatorname{ord}(\eta)=p-1$. Since $\ker\eta$ is
trivial, the corresponding field $L$ is the ray class field of~$K$
mod~$\gp$; that is, $L=H\cdot \Q(\zeta_p)$ (cf.
Propsition~\ref{rayclassfield}). The cyclic group $\Gal (L/H)$ has
order $k:=(p-1)/2$. Also, let $\cE=\End_{H[\Gal(L/H)]}(L)$ be the
$H$-algebra of $\Gal(L/H)$-equivariant endomorphisms. After fixing a
generator $\tau$ of $\Gal(L/H)$, consider~$\Theta$ as in
Lemma~\ref{general}. Finally, let~$\lambda: I(\gp)\rightarrow L^*$
be the one-cocycle in Section~6. To find the elliptic directions in
$S_2(A_f)$ turns out to be equivalent to find the twisted
one-cocycles $\lambda_u(\ga)=\lambda(\ga)\, u/{}^\ga u$ which are
modular. Note that now $\lambda$ is not modular in $[A_f]$.

\bigskip

\begin{proposition}\label{aplication} For all $u\in L^*$, the following conditions are
equivalent:
\begin{itemize}
\item[(i)] the one-cocycle
$\displaystyle{\lambda_u(\ga)=\lambda(\ga)\frac{u}{{}^\ga u}}$ is
modular;
\item[(ii)]
$\displaystyle{u=\Theta\left(\frac{X^k-1}{\Phi_{k}(X)}\right)(v)}$,
for some $\displaystyle{v\not\in \ker
\Theta\left(\frac{X^{k}-1}{\Phi_{k}(X)}\right)}$.
\end{itemize}
In particular, for
$u=\displaystyle{\Theta\left(\frac{X^{k}-1}{\Phi_{k}(X)}\right)(\zeta_p)}$
the one-cocycle $\lambda_u$ is modular. Here, $\Phi_{k}(X)$
denotes the $k$-th cyclotomic polynomial.
\end{proposition}

\noindent {\bf Proof.}  The  values  $u\in L^*$ for which
$\lambda_u$ is modular are the eigenvectors of the $K$-linear map
\begin{equation}
{\operatorname{pr}}(u)=\sum_{\ga\in\Gal(L/K)} {}^{\ga^{-1}}\lambda
(\ga) \left(\sum_{\sigma\in\Phi} \frac{1}{{}^\sigma\psi
(\ga)}\right ){}^{\ga^{-1}} u
\end{equation}
with eigenvalue equal to~$[L: K]$. Also, by
Proposition~\ref{projector}, we know that~${\operatorname{pr}}/[L:
K]$ is a projector, ${\operatorname{pr}}$ diagonalizes, and its
characteristic polynomial is
$$([L: K]-X)^{[E: K]} X^{ [L: K]-[E: K]}=
\left( ([L: K]-X)^{\varphi(k)}X^{k-\varphi(k)}\right)^{[H:
K]}\,.$$ By part (i) of  Lemma~\ref{sumzero}, we can rewrite
$$ {\operatorname{pr}}(u)=\sum_{\ga\in\Gal(L/H)} {}^{\ga^{-1}}\lambda (\ga)
\left(\sum_{\sigma\in\Phi} \frac{1}{{}^\sigma\psi (\ga)}\right ){}^{\ga^{-1}} u\,.$$

Let $g\in \Z$ be a primitive root of $(\Z/ p\,\Z)^*$ such that
$\eta(g)=\zeta$, where~$\zeta=e^{\frac{\pi\,i}{k}}$. Since the set
of principal ideals $\{\ga_j=g^{2 j}\cO_K\colon 1\leq j\leq k \}$
is a set of representatives of $\Gal(L/H)$ and $\lambda (g^{2
j}\cO_K)=g^{2 j}$, we have
$$ G(u):=\frac{{\operatorname{pr}}(u)}{[H: K]}=\frac{1}{[H: K]}\sum_{j=1}^{k}  \left(\sum_{\sigma\in\Phi} {}^\sigma \zeta^{-2 j}\right ){}^{\ga_j^{-1}} u=
\sum_{j=1}^{k}  \Tr_{\Q(\zeta)/\Q} (\zeta^{-2 j})
{}^{\ga_j^{-1}}u\,.$$ Hence,~$G$ belongs to~$\cE$ and its
characteristic polynomial has roots $0$ and $k$ with
multiplicities $k-\varphi(k)$ and $\varphi(k)$, respectively.

Now, we fix the generator $\tau=g^{-2}\cO_K$ of $\Gal(L/H)$ and
apply Lemma~\ref{general} to the endomorphism $G-k
\operatorname{Id}\in \cE$. It follows that the set
$$ \cZ=\{ \zeta'\in \mu_{k} \colon \sum_{j=1}^{k}  \Tr_{\Q(\zeta)/\Q} (\zeta^{-2 j})(\zeta')^{2 j}-k=0\} $$
has cardinality $|\cZ|=\varphi(k)$. Letting $\zeta_k=\zeta^2$, we
claim that
$$\cZ=\{ \zeta_k^{j}\colon  1\leq j < k \,, \gcd(j,k)=1\}\,.$$
Since $\Gal (\Q(\zeta)/\Q)$ acts transitively on $\cZ$ and
$|\cZ|=\varphi(k)$, it suffices to prove that $\zeta_k\in\cZ$.
Indeed, one checks:
$$
\sum_{j=1}^k \left( \sum_{i\in(\Z/k\Z)^*}\zeta_k^{-j\,i}
\right)\zeta_k^j=\sum_{j=1}^k \left(
\sum_{i\in(\Z/k\Z)^*}\zeta_k^{(1-i)\,j} \right)=\sum_{j=1}^k
\left( \sum_{i\in(\Z/k\Z)^*}\zeta_k^{i\,j} \right)=k\,.
$$
Then, from Lemma~\ref{general}, we obtain
$$\{ u\in L\colon  {\operatorname{pr}}(u)=[L: K] \,u \}=
\{ u =\Theta \left(\frac{X^{k}-1}{\Phi_{k}(X)}\right)(v)\colon
v\in L\}\,.
$$
Note that the image of $\Theta
\left(\frac{X^{k}-1}{\Phi_{k}(X)}\right)$ is independent of the
choice of the generator $\tau$ in $\Gal(L/H)$. It can be easily
checked that  $\Theta ((X^{k}-1)/\Phi_{k}(X))$ vanishes on $H$,
which implies that $\Theta ((X^{k}-1)/\Phi_{k}(X))(\zeta_p)$ is
non-zero since the class of the polynomial $(X^{k}-1)/\Phi_{k}(X)$
in $L[X]/(X^k-1)$ is non-zero. \hfill $\Box$ \vskip 0.2 cm

\bigskip

Example: take $p=7$, so that $K=\Q(\sqrt{-7})$ has class number
one. Let~$\psi$ in~$\mathcal X$ with eta-character satisfying
$\eta(3)=e^{2\pi i/6}$. Its corresponding newform~$f=\sum
\psi((a)) q^{\Norm(a)} \in S_2(\Gamma_1(49))$ has
Nebentypus~$\varepsilon$ of order~$3$; note that
$\psi((a))=a\eta(a)$ for all~$a\in \cO_K$. The one-cocyle
$\lambda$ satisfies $\lambda((a))=a$ with the unique choice of
sign for $a$ such that the symbol $(a/\sqrt{-7})=1$. This
one-cocycle is not modular for $\psi$ (in fact, it is modular for
the Hecke character in~$\mathcal X$ with eta-character of
order~$2$ in which case the (unique) elliptic direction coincides
with the rational newform in~$S_2(\Gamma_0(49))$ giving rise to
the elliptic curve $49$A$1$ in Cremona's notation.) Thus, we need
to twist~$\lambda$ by a coboundary in order to get a modular
one-cocycle. According to Proposition~\ref{aplication}, we can
take, for instance, $u=\Theta(X-1)(\zeta_7)=\zeta_7^2-\zeta_7$ and
the cuspidal form $g_{u}=\sum {}^{\ga^{-1}} \lambda_u(\ga)
q^{\Norm(\ga)} = \sum \lambda((a)) {}^{(a^{2})} u / u \,
q^{\Norm(a)} \in S_2(\Gamma_1(49))$ is an elliptic direction
of~$A_f$. A computer calculation shows the lattice $\Lambda$ for
the corresponding elliptic optimal quotient from
$\operatorname{Jac}(X_{\Gamma_\varepsilon})$ satisfies:
$c_4(\Lambda)=c_4(A(7)) u^4$, and $c_6(\Lambda)=c_6(A(7))u^6$.

\bibliographystyle{plain}
\bibliography{MAVCM4}

\end{document}